\documentclass[11pt]{article}
\usepackage{amsfonts,amsmath,amssymb,amsthm,graphicx,float,hyperref,comment,color,subcaption,float,tikz,appendix,url,tabularx,ulem,mathrsfs,url}

\usepackage{bm}
\usepackage[ruled, lined, linesnumbered]{algorithm2e}
\newcommand\revone[1]{{\color{black} #1}}
\newcommand\revonesecondgo[1]{{\color{black} #1}}
\newcommand\revtwo[1]{{\color{black} #1}}
\newcommand\change[1]{{\color{black} #1}}
\newtheorem{definition}{Definition}[section]

\newtheorem{lemma}{Lemma}[section]
\newtheorem{problem}{Problem}[section]

\SetKwComment{Comment}{$\triangleright$ }{}

\newcommand{\ab}{{\mathbf a}}
\newcommand{\bb}{{\mathbf b}}
\newcommand{\cb}{{\mathbf c}}
\newcommand{\db}{{\mathbf d}}
\newcommand{\eb}{{\mathbf e}}

\newcommand{\nb}{{\mathbf n}}

\newcommand{\yb}{{\mathbf y}}

\newcommand{\ub}{{\mathbf u}}
\newcommand{\vb}{{\mathbf v}}

\newcommand{\Wb}{{\mathbf W}}

\newcommand{\xb}{{\mathbf x}}
\newcommand{\Ab}{{\mathbf A}}

\newcommand{\Gb}{{\mathbf G}}

\newcommand{\gamb}{\pmb{\gamma}}

\newcommand{\Cbb}{\mathbb{C}}

\newcommand{\Rbb}{\mathbb{R}}

\newcommand{\Umf}{\mathbb{U}}
\newcommand{\Vmf}{\mathbb{V}}

\newcommand{\nv}{\Vec{n}}

\newcommand{\Tc}{\mathcal{T}}
\newcommand{\Uc}{\mathcal{U}}
\newcommand{\Vc}{\mathcal{V}}
\newcommand{\Sc}{\mathcal{S}}

\newcommand{\Gc}{\mathcal{G}}

\newcommand{\Xc}{\mathcal{X}}
\newcommand{\Wc}{\mathcal{W}}

\newcommand{\Wt}{\widetilde{W}}
\newcommand{\At}{\widetilde{A}}
\newcommand{\Bt}{\widetilde{B}}

\newcommand{\epst}{\widetilde{\epsilon}}

\newcommand{\Abt}{\widetilde{\mathbf{A}}}
\newcommand{\Bbt}{\widetilde{\mathbf{B}}}

\newcommand{\kz}{k_{z_1}}
\newcommand{\kzz}{k_{z_2}}
\newcommand{\kzzz}{k_{z_3}}
\newcommand{\kzt}{\tilde{k}_{z_1}}
\newcommand{\kzzt}{\tilde{k}_{z_2}}
\newcommand{\kzzzt}{\tilde{k}_{z_3}}
\newcommand{\ky}{k_{y_1}}
\newcommand{\kyy}{k_{y_2}}
\newcommand{\kyyy}{k_{y_3}}
\newcommand{\kyt}{\tilde{k}_{y_1}}
\newcommand{\kyyt}{\tilde{k}_{y_2}}
\newcommand{\kyyyt}{\tilde{k}_{y_3}}

\newcommand{\Aa}{\operatorname{AA}}

\newcommand{\LRAAM}{\operatorname{LRAA}(M)}
\newcommand{\PPW}{\operatorname{PPW}}

\usepackage[top=1in,bottom=1in,left=1in,right=1in]{geometry}

\title{LR-WaveHoltz: A Low-Rank Helmholtz Solver}

\author{Andreas Granath\thanks{Department of Mathematics and Mathematical statistics, Ume\r{a} University, Sweden, \href{mailto:andreas.granath@umu.se}{andreas.granath@umu.se}}\and Daniel Appel{\"o}\thanks{Department of Mathematics, Virginia Tech, USA, \href{mailto:appelo@vt.edu}{appelo@vt.edu}}\and Siyang Wang\thanks{Department of Mathematics and Mathematical Statistics, Ume\r{a} University, Sweden, \href{mailto:siyang.wang@umu.se}{siyang.wang@umu.se}}}

\begin{document}

\maketitle
\begin{abstract}
We propose a low-rank method for solving the Helmholtz equation. Our approach is based on the WaveHoltz method, which computes Helmholtz solutions by applying a time-domain filter to the solution of a related wave equation. The wave equation is discretized by high-order multiblock summation-by-parts finite differences. \revtwo{ In two dimensions we seek to compress the solution in matrix form, and in three dimensions using tensor trains}. To control rank growth we use step-truncation during time stepping and a low-rank Anderson acceleration for the WaveHoltz fixed point iteration. We have carried out extensive numerical experiments demonstrating the convergence and efficacy of the iterative scheme for free- and half-space problems in two and three dimensions with constant and piecewise constant wave speeds. 
\end{abstract}

\section{Introduction}
This paper introduces a low-rank method for solving the Helmholtz equation in two and three dimensions. The low-rank techniques we use rely on the use of structured multiblock meshes and the method is therefore best suited to problems posed on domains with low geometric complexity. One area where such problems are commonly found is underwater acoustics. Underwater acoustics is an important field of study due to its wide range of engineering applications, including sonar-based underwater exploration, seismic surveying, and long-distance communication. A key property of acoustic waves that enables these applications is their ability to travel long distances with minimal attenuation. This broad applicability has driven the need for efficient and accurate models of underwater wave propagation. Classical methods for obtaining approximate solutions include the ray tracing method \cite{Hovem2011}, the normal mode method \cite{Pekeris1948}, and parabolic equation methods \cite{Tappert2005}. A common feature of these methods is that they \revtwo{seek to find approximate solutions to the wave propagation problem without solving the Helmholtz equation or the time domain wave equation. In this work we instead try to work directly with compressed representations of the solution to the Helmholtz equation, with the goal to not approximate the physics but still gain efficiency from the compressed format.} Here, \revtwo{we} use the Helmholtz equation
\begin{equation}
    \Delta u(\xb)+\omega^2u(\xb)=f(\xb)\quad\text{ in }\Omega\subset\Rbb^d
    \label{eq:IntroHelmholtz}
\end{equation}
as the underlying PDE, modeling the acoustic pressure $u(\xb)$ of a signal generated by a source $f(\xb)$ with frequency $\omega$ in some spatial domain $\Omega\subset\Rbb^d$ for $d=2,3$. An alternative approach is to solve the Helmholtz equation directly. This is known to be difficult for high frequencies due to the indefinite nature of the equation. In particular, the system of equations obtained by discretizing \eqref{eq:IntroHelmholtz} becomes indefinite \cite{Ernst2011}. This  also limits the effectiveness of iterative solvers: for instance, the conjugate gradient method fails entirely, while GMRES provides little to no acceleration \cite{Erlangga2008}. A further challenge arises from classical dispersion analysis \cite{Kreiss1972,Hagstrom2012}, which shows that discretizing a wave problem with frequency $\omega$ using a finite difference method of order $2p$ and a prescribed accuracy $\epsilon$ requires a number of points per wavelength (PPW) on the order of $\PPW \sim(\omega / \epsilon)^{1/2p}$. It is therefore more efficient to use higher order methods.

The WaveHoltz method was introduced in \cite{Appelo2020} as an alternative approach for solving \eqref{eq:IntroHelmholtz}. Specifically, the problem is solved in the time domain over one period, $T=\frac{2\pi}{\omega}$ and the result is then filtered. The filtered solution serves as the initial data for the next iteration. When the WaveHoltz system is discretized, the resulting matrix is positive definite, unlike the indefinite system obtained from a direct frequency domain discretization. This property makes Krylov subspace methods suitable for accelerating the iteration. In \cite{Appelo2020}, both CG and GMRES are applied in the context of the acoustic wave equation; CG is used for the elastic wave equation in \cite{Appelo2020elastic}, and GMRES is employed for Maxwell's equations in \cite{Peng2022}.

The finite difference method has seen a lot of success in the context of solving higher-order wave problems. However, a known limitation is the difficulty of obtaining stability estimates for initial-boundary value problems. To introduce additional structure, finite difference operators satisfying the summation-by-parts (SBP) property were developed \cite{Kreiss1974}. These operators satisfy a discrete analogue of the integration-by-parts identity, enabling the derivation of energy estimates analogous to those in the continuous setting. The SBP operators themselves do not enforce boundary conditions; instead, they are often combined with the simultaneous approximation term (SAT) method, which weakly enforces boundary conditions through penalization \cite{Carpenter1994}. The SBP-SAT framework has been extensively studied and applied, including to the two-dimensional acoustic wave equation in second-order form, both in continuous media \cite{Mattsson2009} and in discontinuous media with material interfaces \cite{Mattsson2008}.

To further decrease the complexity of representing the solution, there has been a large interest in low-rank methods for discretizing differential equations. The underlying idea is to utilize, when present,  inherent low-rank structures in the solution to reduce the cost of storage and computation \cite{Bachmayr2023}. \revtwo{In two dimensions on structured block meshes, we represent the solution as a compressed matrix using the dyadic SVD form, and in three dimensions we use the \textit{tensor train} (TT) format \cite{Oseledets2011}.} 

\revtwo{Suppose in $d$ dimensions the solution is discretized using $n^d$ gridpoints. The goal of this paper is to design methods whose computational complexity scale like $\mathcal{O}(d n r^k)$, where $k$ is between 2 and 4, for solutions that have a matrix or tensor rank $r$. To do this we use a low rank discretization of the wave equation. The timestepping method is of the step-truncation type, where the rank is allowed to grow during the time step and then truncated using SVD with a given tolerance, see for example \cite{Dektor2021}. An alternative could be to use the \textit{dynamic low-rank} approximation method \cite{Koch2007dynamical}, where the rank of the solution is constant. The  
WaveHoltz fixed-point iteration is then either iterated directly or accelerated using a low-rank formulation of Anderson acceleration (LRAA) \cite{Appelo2025}.}

Naturally, a low-rank method is only computationally efficient if the solution itself exhibits low-rank structure. While this is often the case for coercive elliptic problems \cite{LowRankFE}, one should not generally expect the solution to the Helmholtz equation to be low-rank. In fact, the work by Engquist and Zhao \cite{no_low_rank_Engq} demonstrates that, in some cases, the separation rank grows rapidly with frequency. However, there are scenarios where low-rank structures in Helmholtz solutions can be exploited. A prominent example is the fast directional algorithm for Helmholtz scattering problems introduced by Engquist and Ying \cite{parabolic_sep}. This algorithm leverages the low-rank separability of the Green's function when considering interactions between a ball of radius $r$
and a distant region located at a distance 
 $r^2$, confined within a cone of opening angle $1/r$. In this work, we argue that similar low-rank structures can be exploited in volume discretizations for underwater acoustics problems involving a single point source.

Adding to the discussion above, we empirically investigate the potential for a fast volume discretization method. We present a low-rank framework for solving the Helmholtz equation in two- and three-dimensional Cartesian domains using the WaveHoltz method, which we refer to as the \textit{LR-WaveHoltz} (LRWH). The LRWH method builds on the WaveHoltz framework and incorporates a multi-block SBP-SAT solver, with all unknowns stored in low-rank form, using the SVD in two dimensions and the TT format in three dimensions. The multi-block structure is essential, as the point source is not well separated from all regions of the domain. In blocks that satisfy the "parabolic separation" condition of \cite{parabolic_sep}, we expect the solution rank to remain low; for blocks near the source, however, the rank is likely to be high. The method proposed here should therefore be regarded as a proof of concept. A more efficient implementation would likely involve hybridizing traditional solvers near the source with low-rank solvers in regions farther away, an approach that loosely corresponds to the separation of near- and far-field terms in multipole methods.

The rest of the paper is organized as follows: in Section~\ref{sec:waveholtz} we introduce the model problem and WaveHoltz method, as well as the Anderson Acceleration algorithm. In Section~\ref{sec:fullranksection} we introduce the SBP method and full rank semidiscretization of the model problem. Then, in Section~\ref{sec:LRSBPsection} we introduce the notion of step truncation, as well as the low-rank wave solver in two and three dimensions. We then extend this to the LRWH method in Section~\ref{sec:LRWHsection} and end the paper by demonstrating some numerical examples in Section~\ref{sec:numericalexamples} and conclude in Section ~\ref{sec:conclusion}.

\section{Model problem and the WaveHoltz method}
In this section we begin by introducing the model problem, defining the physical setting which we are interested in and the associated notation. We then present the WaveHoltz method that will be used to solve the model problem. Finally, we end the section by introducing the Anderson acceleration algorithm which will be used in conjunction with the WaveHoltz method.
\label{sec:waveholtz}
\subsection{The model problem}
We consider a smooth open domain $\Omega\subset\Rbb^d$ which is assumed to either be fully submerged under water or have its top boundary as the water surface. The problem is driven by a forcing $f(\xb)$, has a potentially discontinuously varying wave speed $c$ and a damping $\kappa(\xb)$. We are interested in the distribution of the acoustic pressure $u(\xb)$ generated by a forcing oscillating at a given frequency $\omega$, which is modeled by the Helmholtz equation \cite{Jenssen2011}
\begin{alignat}{2}
    &\nabla\cdot(c^2\nabla u(\xb))+\omega^2u(\xb)-i\omega\kappa(\xb)u=f(\xb)&\text{ in }\Omega
    \label{eq:Helmholtz},\\
    &  ia\omega u(\xb)+bc^2\nabla u(\xb)\cdot\nv=0, &\text{ on }\partial\Omega,
    \label{eq:Boundarycondition}
\end{alignat}
where $a,b$ are constants satisfying $a^2+b^2=1$. 
The boundary conditions \eqref{eq:Boundarycondition} model the lowest order nonreflecting boundary conditions when $a,b=\frac{1}{\sqrt{2}}$. 

 Extensions of the model could take higher order nonreflecting boundary conditions \cite{Hagstrom1999} and acoustic-elastic interaction with the sea floor \cite{Appelo2019} into account.
\subsection{The WaveHoltz method}

 We briefly summarize the WaveHoltz method developed in \cite{Appelo2020}. Introducing the $T=\frac{2\pi}{\omega}$-periodic function $w=w(\xb,t)=e^{i\omega t}u(\xb)$, the associated wave problem of \eqref{eq:Helmholtz}-\eqref{eq:Boundarycondition} is then given by

\begin{alignat}{2}
    &w_{tt}(\xb,t)+\kappa(\xb)w_t(\xb,t)=\nabla\cdot(c^2\nabla w(\xb,t))-f(\xb)\cos(\omega t)\quad&&\text{ for }\quad (\xb,t)\in\Omega\times (0,T],\label{eq:waveproblem} \\
    &aw_t(\xb,t)+bc^2\nabla w(\xb,t)\cdot\nv=0\quad &&\text{ for }\quad (\xb,t) \in\partial\Omega\times (0,T],\nonumber\\
    & w(\xb,0)=v_0(\xb), \quad w_t(\xb,0)=v_1(\xb)\quad &&\text{ for }\quad \xb\in\Omega,\nonumber    
    \label{eq:waveproblem}
\end{alignat}
where $v_0(\xb)=u(\xb)$ and $v_1(\xb)=i\omega u(\xb)$. 
Assuming that the initial data of the wave problem satisfies $v_0\in H^1(\Omega)$ and $v_1\in L^2(\Omega)$, the WaveHoltz (WH) operator $\Pi$ acting on the data is defined as
\begin{equation}
    \Pi\begin{bmatrix}
        v_0(\xb) \\ v_1(\xb)
    \end{bmatrix}=\frac{2}{T}\int_0^T\bigg(\cos(\omega t)-\frac{1}{4}\bigg)\begin{bmatrix}
        w(\xb,t)\\ w_t(\xb,t)
    \end{bmatrix}dt,\quad T=\frac{2\pi}{\omega}.
    \label{eq:operator}
\end{equation}
The kernel in $\Pi$ is chosen to damp the frequencies away from $\omega$, while exactly preserving the $\omega$ mode. Applying this operator can therefore be viewed as a filtering of frequencies around $\omega$. Evaluating the integral in \eqref{eq:operator} results in $[u(x), i\omega u(\xb)]^T$, hence the Helmholtz solution $u(\xb)$ solves the fixed point problem $\Pi [u(\xb), i\omega u(\xb)]^T=[u(x), i\omega u(x)]^T$.  The WaveHoltz method is then the fixed point iteration
\begin{equation}
    \begin{bmatrix}
        v(\xb) \\
        v(\xb)'
    \end{bmatrix}^{(n+1)}=\Pi\begin{bmatrix}
        v(\xb)\\
        v(\xb)'
    \end{bmatrix}^{(n)}\quad\text{ with }\quad \begin{bmatrix}
        v(\xb)\\ v(\xb)'
    \end{bmatrix}^{(0)}=\begin{bmatrix}
        0 \\ 0
    \end{bmatrix}.
    \label{eq:zeroiteration}
\end{equation}
The convergence of the WaveHoltz method has been studied in \cite{Appelo2020,ROTEM2026} where it was shown that under natural assumptions the iteration converges to the Helmholtz solution, i.e.  $u=\lim_{n\rightarrow\infty} v^{(n)}$. \revone{It is expected that the number of iterations to converge grows like $\mathcal{O}(\omega)$, discussed in the context of impedance boundary conditions in \cite{ROTEM2026}}.

\begin{algorithm}[ht!]
\caption{Anderson Acceleration}
\label{algorithm:AA}
\KwIn{Initial condition $X^0\in\Xc$, memory parameter $m$ and error tolerance $\epsilon^\star$}
\KwOut{ $X^k\in\Xc$ approximately solving $X^k=G(X^k)$ within tolerance $\epsilon^\star$}
\For{$k=1,2,...$\text{ until } $\|X^{k+1}-X^k\|<\epsilon^\star$}{
$m_k=\min(m,k)$\\
Set $D_k=[\Delta F_{k-m_k},...,\Delta F_{k-1}]^T\in\Rbb^{m_k}$, where $\Delta F_{i}=F_{i+1}-F_i$ for $F_i=G(X^i)-X^i$\\
Solve for $\pmb{\gamma}^{(k)}=[\gamma_0^{(k)},\gamma_1^{(k)},...,\gamma_{m_k-1}^{(k)}]^T\in\Rbb^{m_k}$ satisfying
\begin{equation}
    \pmb{\gamma}^{(k)}=\operatorname{argmin}_{\pmb{u}\in\Rbb^{m_k}}(\|D_k\pmb{u}-F_k\|)
    \label{eq:AAweights}
\end{equation}\\
Update $X_k$ according to
\begin{equation}
    X^{k+1}=\gamma^{(k)}_0G(X^{k-m_k})+(1-\gamma^{(k)}_{m_k-1})G(X^k)+\sum_{i=0}^{m_k-1}(\gamma^{(k)}_{i+1}-\gamma^{(k)}_i)G(X^{k-m_k+i})
    \label{AA:Gsum}
\end{equation}\\
Calculate $\|X^{k+1}-X^k\|$
}
\end{algorithm}
\subsection{Anderson acceleration}
As demonstrated in \cite{Appelo2020} and \cite{ROTEM2026}, the WaveHoltz iterations can be accelerated using a Krylov subspace method. One such method suitable for our low-rank formulation is \textit{Anderson acceleration} (AA). Following the notation in \cite{Yang2022}, let $G:\Xc\rightarrow\Xc$ for some Hilbert space $\Xc\subset L^2(\Omega)$, $X\in\Xc$ and consider the fixed point problem $G(X)=X$. This is commonly solved using Picard iteration; given some initial condition $X^0\in\Xc$, update the iterates using
$$X^{k+1}=G(X^k)\quad\text{ for }k=1,2,3... $$
until convergence. Picard iteration only uses information from the previous iterate and converges conditionally. Moreover, the convergence is typically linear.  To accelerate the convergence, AA instead uses a linear combination of a given number of previous iterates. More precisely, given a \textit{memory parameter} or \textit{window size} $m$, the iterate $X^k$ of the $\operatorname{AA}(m)$ method will use the previous $\min(m,k)+1$ iterates, outlined in Algorithm~\ref{algorithm:AA}. From this definition we can infer that $\Aa (0)$ corresponds to Picard iteration.

Note that $\|\cdot\|$ denotes a general norm in \cite{Yang2022}, but we resort to the Frobenius norm. The naive way of solving for the coefficients $\pmb{\gamma}^{(k)}$ when $X^k\in\Rbb^{n\times n}$ in \eqref{eq:AAweights} requires solving a large system of linear equations of size $\mathcal{O}(mn^2)$ at each iteration. We reformulate the minimization problem in terms of a smaller system of the size $\mathcal{O}(m^2)$ in Problem~\ref{Lemma:AAweights}.
\begin{problem}
\label{Lemma:AAweights}
    Let $\|\cdot\|$ denote the Frobenius norm, $\{\Delta F_{k-j}\}_{j=1}^{m}$ with $\Delta F_l\in\Rbb^{m}$ be given as well as $F_k\in\Rbb^{m}$, then the vector $\pmb{\gamma}^{(k)}$ satisfying
    $$\pmb{\gamma}^{(k)}=\operatorname{argmin}_{\ub\in\Rbb^m}\|D_k\ub-F_k\|,$$
    can be obtained by solving $A\pmb{\gamma}^{(k)}=b$ with
    $$A_{ij}=\langle\Delta F_i,\Delta F_j\rangle,\quad b_i=\langle\Delta F_i,F_k\rangle.$$
\end{problem}

\section{Full rank wave solver in matrix form}
\label{sec:fullranksection}
We begin the section by providing a background on the numerical method used for the spatial discretization of \eqref{eq:Helmholtz}-\eqref{eq:Boundarycondition}, after which we introduce the full rank semidiscretization of the problem, which is then generalized to the low-rank setting in Section~\ref{sec:LRmethods}. 
\subsection{SBP operators}
\label{subsec:SBP}
Let $I=[0,1]$ denote a one-dimensional interval discretized uniformly into $n$ gridpoints $\{x_i\}$ with grid spacing $h$ such that $x_i=(i-1)h$ and $h=\frac{1}{n-1}$. Moreover, let $u$ and $v$ be functions defined on $I$ with sufficient regularity, denoting the evaluation at a point $x_i$ as $u_i=u(x_i),v_i=v(x_i)$ and grid evaluations as $\ub=[u_1,...,u_n]^T,\vb=[v_1,...,v_n]^T$. We introduce the notion of narrow-stencil SBP operators $D_1\approx\partial_x$ and $D_2\approx\partial_x^2$ having order of accuracy $\mathcal{O}(h^{2p})$ in the interior and $\mathcal{O}(h^p)$ at a few points near the boundary. The decreased order of accuracy is due to the one-sided stencil and we use the name convention that the operators are of order $2p$. We can then define the SBP property of the operators $D_1$ and $D_2$ \cite{Mattsson2004}.
\begin{definition}
     The operators $D_1$ and $D_2$  are said to satisfy the summation-by-parts (SBP) property if
    $$HD_1+(HD_1)^T=\eb_n\eb_n^T-\eb_1\eb_1^T,$$
    $$D_2=H^{-1}(-A+(\eb_n\eb_n^T-\eb_1\eb_1^T)S),$$
    where $\eb_1=[1,\cdot\cdot\cdot,0]^T,\eb_n=[0,\cdot\cdot\cdot,1]^T$, $H$ is symmetric positive definite, $A$ is symmetric positive semidefinite and $S$ is an approximation of the derivative at the boundary.
    \label{def:SBPproperty}
\end{definition}
The matrix $H$ defines a quadrature such that $\ub^TH\vb\approx\int_Iu(x)v(x)dx$ and an associated inner product \cite{Hicken2013}. It should also be noted that the boundary derivative approximations have an improved order of accuracy $\mathcal{O}(h^{p+1})$ compared with $D_1$. 

\subsection{The full rank formulation in 2D}
To discretize in space, we resort to the method outlined in \cite{Mattsson2008}, representing the numerical solution in each block as a matrix rather than a vector. For ease of presentation we restrict ourselves to a two-block domain with aligning grids, but the framework readily extends to a multiblock setting. Let $\Omega\subset\Rbb^2$ with be a Cartesian domain, with coordinates $(x_1,x_2)\in\Omega$, partitioned into two subdomains $\Omega_u$ and $\Omega_v$ having a vertical interface $\Gamma=\partial\Omega_u\cap\partial\Omega_v$ at $x_1=0$ and constant, but not necessary equal, associated wave speeds $c_u$ and $c_v$. Furthermore, we denote the positive outward unit normal along $\partial\Omega$ as $\nb$. The continuous problem in $\Omega$ can then be formulated as
\begin{alignat*}{2}
    &u_{tt}(\xb,t)=\nabla\cdot(c_u^2\nabla u(\xb,t))-f(\xb)\cos(\omega t),&&\text{ for }\quad (\xb,t)\in\Omega_u\times(0,T),\\
    &v_{tt}(\xb,t)=\nabla\cdot(c_v^2\nabla v(\xb,t))-f(\xb)\cos(\omega t),&&\text{ for }\quad(\xb,t)\in\Omega_v\times(0,T),\\
    &u_t(\xb,t)+c_u\nabla u(\xb,t)\cdot\nb=0,&&\text{ on }\partial\Omega_u\backslash\Gamma,\\
    &v_t(\xb,t)+c_v\nabla v(\xb,t)\cdot\nb=0,&&\text{ on }\partial\Omega_v\backslash\Gamma,\\
    &u(\xb,0)=0,\quad u_t(\xb,0)=0, &&\text{ for }x\in\Omega_u,\\
    &v(\xb,0)=0,\quad v_t(\xb,0)=0, &&\text{ for }x\in\Omega_v,
\end{alignat*}
where letters in the subscripts denote the domains. The solution and its normal flux across the interface must be continuous. Letting $\nb_\Gamma$ denote the unit normal $\nb_\Gamma=[1,0]^T$ at $\Gamma$, we formulate these conditions in terms of $u$ and $v$ as 
\begin{alignat}{2}
&u=v&&\text{ at } x_1=0\label{eq:interfacecond1},\\
&c_u^2\nabla u\cdot\nb_\Gamma-c_v^2\nabla v\cdot\nb_\Gamma=0&&\text{ at } x_1=0\label{eq:interfacecond2}.
\end{alignat}

We now turn to the SBP-SAT semidiscretization of the problem. Assume that the domains are Cartesian blocks discretized by $n$ degrees of freedom $\xb,\yb\in\Rbb^n$ in each direction. Then, we represent the time dependent grid evaluation of $u,v$ in each block as matrices $\Umf(t),\Vmf(t)\in\Rbb^{n\times n}$ ordered such that $\Umf_{ij}(t)=u(x_i,y_j,t)$ and $\Vmf_{ij}(t)=v(x_i,y_j,t)$. Similarly, let $F_u,F_v\in\Rbb^{n\times n}$ denote the grid evaluations of the forcing $f(x,y)$ in the corresponding block. We use the difference operators defined in the previous section to approximate spatial derivatives as $\partial_x^2u(t)\approx D_2\Umf(t)$ and $\partial_y^2u(t)\approx \Umf(t)D_2^T$. The semidiscretization is then given by
\begin{align}
    &\Umf_{tt}=c_u^2D_2\Umf+c_u^2\Umf D_2^T+\sum_{i=1}^4(-1)^i\Sc^u_i(\Umf,\Vmf)-F_u\cos(\omega t),
    \label{eq:fullrankdiscU}\\
&\Vmf_{tt}=c_v^2D_2\Vmf+c_v^2\Vmf D_2^T+\sum_{i=1}^4(-1)^i\Sc^v_i(\Umf,\Vmf)-F_v\cos(\omega t),
\label{eq:fullrankdiscV}
\end{align}
where the first two terms approximate the Laplacian, the sum over SAT terms $\Sc^u_i(\Umf,\Vmf),\Sc^v_i(\Umf,\Vmf)$ enforces interface or boundary conditions and the last term is the periodic forcing function. The indices of the SAT terms are ordered to represent different boundaries of the respective block, where the order is given as west, north, south and east. We present two of the SAT terms in the domain $\Omega_u$ as the remaining terms are constructed analogously. At the northern boundary, the nonreflecting condition $u_t+c_u\nabla u\cdot\nb=0$ with $\nb=[0,1]^T$ is imposed weakly by the term $\Sc_2^u(\Umf,\Vmf)$ defined by  
\begin{equation}
    \mathcal{S}_2^u(\Umf,\Vmf)=-(c_u\Umf S^T+\Umf_t)(c_uH^{-1}\eb_n\eb_n^T)^T,
\end{equation}
where we recall the definition of $H,S$ and $\eb_n$ from section~\ref{subsec:SBP}. Again, note that the operators act from the right on $\Umf$ as the boundary condition is in the vertical direction. Since the grid within the blocks are assumed to align, we impose the interface conditions \eqref{eq:interfacecond1}-\eqref{eq:interfacecond2} weakly along $\Gamma$, following the method outlined in \cite{Mattsson2008}. \revone{ If the mesh size is different in the two blocks, which might be desireable when the wave speed varies between the blocks, we can use the standard SBP projection approach from \cite{Almquist2019}.} In $\Omega_u$, the interface conditions are imposed by the SAT term $\Sc_4(\Umf,\Vmf)$ 
\begin{equation*}
\mathcal{S}_4^u(\Umf,\Vmf)=\frac{c_u^2}{2}H^{-1}S^T\eb_n(\eb_n^T\Umf-\eb_1^T\Vmf)-\frac{1}{2}H^{-1}\eb_n(c_u^2\eb_n^TS\Umf-c_v^2\eb_1^TS\Vmf)-c_{av}\frac{\tau}{h}H^{-1}\eb_n(\eb_n^T\Umf-\eb_1^T\Vmf),
\end{equation*}
where $\tau>0$ denotes a stability parameter independent of wave speeds chosen suitably large depending on the order of the SBP method used, and $c_{av}=\frac{c_u^2+c_v^2}{2}$ \cite{Mattsson2008}. Since the semi discretization is a matrix formulation of the method presented in \cite{Mattsson2008} their stability result also applies here. It therefore follows that \eqref{eq:fullrankdiscU}-\eqref{eq:fullrankdiscV} satisfies an energy estimate when choosing $\tau$ sufficiently large. 

\section{The step-truncated LR-SBP method}
\label{sec:LRSBPsection}
We begin by introducing the notion of explicit step truncation methods for general matrices \cite{Dektor2021}. Then, we construct a two dimensional low-rank semidiscretization of \eqref{eq:waveproblem} using the explicit step truncation. Finally, we end the section by extending the method to three dimensions, representing the solution with tensor trains rather than SVD.
\label{sec:LRmethods}
\subsection{Explicit step-truncation }
\label{subsec:step-truncation}
 Let $W(t)\in\Rbb^{n\times n}$ denote the grid evaluation of a function $w$. We consider a general problem on the form
 \begin{equation}
     \frac{\partial^2 W}{\partial t^2}=G(W,t),
     \label{eq:modelproblem}
 \end{equation}
for some right-hand side $G(\cdot,\cdot)\in\Rbb^{n\times n}$ being linear in its first argument. Denote the discrete solution to \eqref{eq:modelproblem} at time $t^k$ by $W^k$ and let its singular value decomposition (SVD) be given as $W^k=U^kS^k(V^k)^T$, where $U^k\in\Rbb^{n\times r}, S^k\in\Rbb^{r\times r}, V^k\in\Rbb^{n\times r}$ for some $r>0$ denoting the \textit{rank} of $W^k$. We introduce the \textit{truncated SVD operator} $\Tc_\epsilon$ which approximates $W^k$ within a given accuracy $\epsilon$ as $\|W^k-\Tc_\epsilon(W^k)\|<\epsilon$ where $\|\cdot\|$ denotes the Frobenius norm. The truncation tolerance $\epsilon$ is chosen such that the residual at the final iteration does not exceed the local truncation error of the underlying spatial discretization. The notation $\Tc_\epsilon(W^k)$ should be understood as the operator acting on the SVD factors of $W^k$.  

When we discretize \eqref{eq:modelproblem} we will add up several terms at each time step. This addition will introduce rank growth. The problems we are targeting are assumed to have low rank at all times, we therefore need a summation operation $\Tc_\epsilon^{sum}$ mimicking the properties of $\Tc_\epsilon$, summarized in Algorithm~\ref{algorithm: summation} using MATLAB notation. We note that if the number of summands is large, the algorithm can be made more efficient by making it sequential. However, since we only have a few terms in the time discretization, we use the algorithm as presented.

\begin{algorithm}[H]
 \caption{Sum of low-rank matrices} 
 \label{algorithm: summation}
\KwIn{low-rank matrices $W_j$, $j=1,\ldots,N$ in SVD form $U_j S_j (V_j)^T$, $j=1,\ldots,N$, and tolerance $\epsilon$}
\KwOut{truncated SVD of sum, $\Tc_\epsilon^{sum}(\sum_{j=1}^N W_j)=\Uc\Sc\Vc^T$}
Form $U=[U_1,\ldots,U_N], S=\operatorname{diag}(S_1,\ldots,S_N), V=[V_1,\ldots,V_N]$\\
Form pivoted QR decompositions: $[Q_1,R_1,\Pi_1]=\textbf{qr}(U), [Q_2,R_2,\Pi_2]=\textbf{qr}(V)$\\
Compute truncated SVD: $\Tc_\epsilon(R_1\Pi_1 S\Pi_2^TR_2^T)=\Uc\Sc\Vc^T$\\
Form $\Uc\gets Q_1\Uc$, $\Vc\gets Q_2\Vc$\\
Return $[\Uc,\Sc,\Vc]=\Tc_\epsilon^{sum}(\sum_{j=1}^NW_j)$
\end{algorithm}
We resort to second order centered differences to discretize \eqref{eq:modelproblem} in time, resulting in \eqref{eq:FRSylvester} and present the evolution of one time step in Algorithm~\ref{algorithm:LRleapfrog}. The first step of the method consists of a regular leapfrog step for \eqref{eq:modelproblem} with a truncated right hand side, yielding a solution $\Wt^{k+1}$ which may have a larger rank than desired. The second step then consists of truncating the sum in the first step using the summation operator $\Tc_{\epsilon}^{sum}(\cdot)$ described in Algorithm~\ref{algorithm: summation}.  
\begin{algorithm}
    \caption{low-rank time evolution from $t^k$ to $t^{k+1}$}
    \label{algorithm:LRleapfrog}
    \KwIn{numerical solution at time $t^k$ in SVD form $W^k=U^kS^k(V^k)^T$, time step $\Delta t$, truncation tolerance $\epsilon$}
    \KwOut{numerical solution at time $t^{k+1}$ $W^{k+1}=U^{k+1}S^{k+1}(V^{k+1})^T$}
    $\Wt^{k+1}=2W^k-W^{k-1}+(\Delta t)^2\Tc^{sum}_{\epsilon}(G(W^k,t^k))$\\
    $W^{k+1}=\Tc_{\epsilon}^{sum}(\Wt^{k+1})$\\
    Return $W^{k+1}=U^{k+1}S^{k+1}(V^{k+1})^T$
\end{algorithm}
\subsection{The two dimensional SVD representation}
\label{subsec:LRSBP}
\begin{figure}[H]
    
    \begin{subfigure}{0.47\linewidth}
    \centering
          \includegraphics[width=0.75\linewidth]{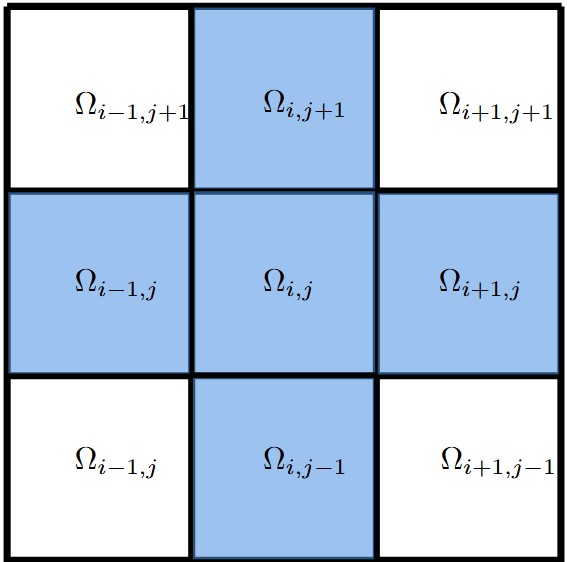}
          \subcaption{}
              \label{fig:domainpartition}
    \end{subfigure}
        \begin{subfigure}{0.6\linewidth}
        \centering
          \includegraphics[width=0.6\linewidth]{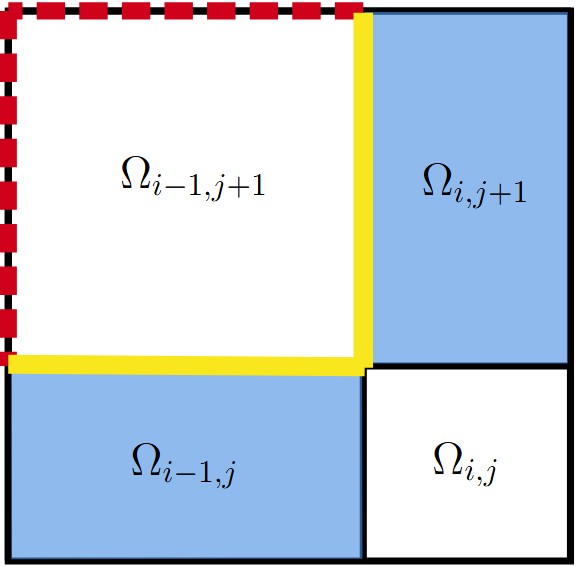}
          \subcaption{}
          \label{fig:cornerblock}
    \end{subfigure}
  
    \caption{(a) An illustration of a multiblock partitioning of $\Omega$, highlighting the multiblocks contributing terms to the discrete Laplacian $L_{ij}$ via interface conditions on $\partial\Omega_{ij}$. (b) An illustration of a corner block, indicating the boundaries (dashed) contributing to the coefficient matrices in the Sylvester equation \eqref{eq:FRSylvester2} via the boundary conditions. We also highlight the internal edges where interface conditions are imposed (solid).}

\end{figure}

 We now generalize the discretization of the domain $\Omega$ used in Section~\ref{subsec:SBP} to a general multiblock setting. For notational convenience, we assume that $\Omega$ is partitioned into $m$ blocks $\Omega_{ij}$ in each direction, where the subscripts from now on do not denote elements, but rather the indices of a block. The blocks are then discretized with $n$ gridpoints and grid size $h$ in each dimension. We denote the grid evaluation of the acoustic pressure $w$ and velocity $w'$ in each block as the matrices $W_{ij}(t),W'_{ij}(t)\in\Rbb^{n\times n}$ for $i,j=1,\ldots,m$. We also need the operator $L_{ij}(W_{ij})\in\Rbb^{n\times n}$ , referred to as the discrete Laplacian, combining $D_2W_{ij}+W_{ij}D_2^T$ with the factors from SAT-terms imposing boundary and interface conditions on $\partial\Omega_{ij}$. An illustration of the contributing edges is shown in Figure~\ref{fig:domainpartition}. The main idea is that given initial data in SVD form, we  construct a scheme that retains the low-rank structure of the solution at all iterations, without ever having to form the full matrix. This can be done in a straightforward manner following the semidiscretization shown in \eqref{eq:fullrankdiscU}-\eqref{eq:fullrankdiscV} where the discrete Laplacian operates on the left or right singular vectors of the numerical solution. We present the solver in Algorithm~\ref{algorithm:LRWaveSolver} for a time interval $[0,T]$ discretized with time step $\Delta t$ and discuss it further below.
\begin{algorithm}[ht!]
 \caption{low-rank wave equation solver} 
 \DontPrintSemicolon
 \label{algorithm:LRWaveSolver}
 \KwIn{Initial acoustic pressure $\{W_{ij}^0\}_{i,j=1}^m$ and velocity $\{{W^0}'_{ij}\}_{i,j=1}^m$ in SVD form, combined Laplacian and SAT operators $\{L_{ij}\}_{i,j=1}^m$ , discretized forcing $\{F_{ij}\}_{i,j=1}^m$ in SVD form, time step $\Delta t$, final time $T$, grid size $h$, truncation tolerance $\epsilon$, frequency $\omega$}
 \KwOut{Approximate solutions $W_{ij}$ to \eqref{eq:waveproblem} at time $T$ in SVD form for each block}
  \For{$ij\in\text{ all multiblocks}$}{
        $W^{-1}_{ij}=\Tc^{sum}_{\epsilon}[W^0_{ij}-\Delta t{W^0}'_{ij}+\frac{(\Delta t)^2}{2}\Tc^{sum}_{\epsilon}(L_{ij}(W^0_{ij}))]$\Comment*[r]{Solution at $-\Delta t$}
  
    }
 \For{$k=0,1,2,3,\ldots$\textnormal{ until } $k\Delta t=T$}{
 $t^k=k\Delta t$\\
    \For{$ij\in\text{ all multiblocks }$}{

            $\widehat{W}_{ij}^{k+1}=\Tc_{\epsilon}^{sum}[L_{ij}(W^k,t^k)+F_{ij}\cos(\omega t^k)]$\Comment*[r]{Sum intermediate terms}
            $W^{k+1}_{ij}=\Tc_{\epsilon}^{sum}[2W^k_{ij}-W^{k-1}_{ij}+(\Delta t)^2\Tc_{\epsilon}^{sum}(\widehat{W}_{ij}^{k+1})]$\Comment*[r]{Evolve from $t^k$ to $t^{k+1}$}
           $W^{k+1}_{ij}=\texttt{LRfADI}(W^{k+1}_{ij},i,j)$\Comment*[r]{If physical corner or corner boundary condition}}        
}
\end{algorithm}

The method follows a standard wave solver using the leapfrog scheme in time, with the difference being the need of applying the truncation operators. We again stress that the notation $W_{ij}^k,W_{ij}'$ and $F_{ij}^k$ should be understood as in SVD form and all operations are done on the factors of the solution and never on the full matrix. The first step of the method consists of the usual construction of the solution at time $-\Delta t$ required for leapfrog method. Then, looping over each block, we first compute $\widehat{W}^{k+1}_{ij}$, consisting of the discrete Laplacian and forcing term. We then add the contribution from the second derivative to form $\Wt^{k+1}_{ij}$ and round with an overall tolerance of $\epsilon$.

We now consider the terms arising from $W'(t)$ in the nonreflecting boundary conditions. Consider the corner block in Figure~\ref{fig:cornerblock}. For notational convenience, let $c$ denote the wave speed in this multiblock. When $W^{k+1}_{ij}$-terms  are combined we obtain
$$W^{k+1}_{ij}+\change{\frac{c\Delta t}{2}}H^{-1}(\eb_0\eb_0^T)W^{k+1}_{ij}-\change{\frac{c\Delta t}{2}}W^{k+1}_{ij}(\eb_n\eb_n^T)H^{-1}=R,$$
where the matrix $R$ depends on $W^k_{ij}$ and $W^{k-1}_{ij}$. A difference compared to vectorized formulations is now that $W^{k+1}_{ij}$ cannot be obtained by a straightforward linear solve. Instead, we have to solve a matrix equation. To this end, we recast the problem into a Sylvester equation $A W^{k+1}_{ij}-W^{k+1}_{ij}B^T=R$, which we then solve using the \textit{factored Alternating Direction Iteration} \cite{Benner2009}. In this context fADI converges in exactly \revone{three} iterations. This method is particularly useful as the solution is immediately expressed in SVD form. A more detailed description of the method is provided in \revone{\ref{subsec:fADI}}. 

\subsection{The three dimensional tensor train representation}
We now outline the low-rank three dimensional discretization of \eqref{eq:waveproblem}. On a multiblock in three dimensions the numerical solution is expressed in the form of an array with elements $\Ab(i_1,i_2,i_3)$ for $i_k=1,2,\ldots,n_k$. As in two dimensions we never represent all of $\Ab$ but rather use a compressed version $\Abt$. Precisely, to represent the compressed grid functions in three dimensions we use the tensor train (TT) format. Here we introduce basic TT concepts, for a more detailed introduction we refer the reader to \cite{Oseledets2011}.  

As in the two dimensional case we use the Frobenius norm. Thus the  inner product $\langle\cdot,\cdot\rangle$ and induced norm $\|\cdot\|$ are
\[
\langle\Abt,\Bbt\rangle=\sum_{i_1,i_2,i_3}\Abt(i_1,i_2,i_3)\Bbt(i_1,i_2,i_3),\quad \|\Abt\|=\sqrt{\langle \Abt,\Abt\rangle},
\]
where $\|\cdot\|$ corresponds to the standard Frobenius norm.  The TT-rank, $(r_0,r_1,r_2,r_3)$, TT-representation of $\Abt$ is then
\[
\Abt(i_1,i_2,i_3)=\sum_{\alpha_0=1}^{r_0}\sum_{\alpha_1=1}^{r_1}\sum_{\alpha_2=1}^{r_2}\sum_{\alpha_3=1}^{r_3}\Gc_1(\alpha_0,i_1,\alpha_1)\Gc_2(\alpha_1,i_2,\alpha_2)\Gc_3(\alpha_2,i_3,\alpha_3), \ \  i_k=1,2,\ldots,n_k,k\in\{1,2,3\}
\]
The first and the last rank is always one, i.e $r_0=r_3=1$. The factor matrices $\Gc_k$ are called \textit{TT-cores} and $\Gc_k\in\Rbb^{r_{k-1}\times n_k\times r_k}$. It is then clear that the cost of storing $\Abt$ is bounded by $3nr^2$, where $r=\max\{r_1,r_2\}$ and $n=\max\{n_1,n_2,n_3\}$, which is smaller than the cost $n^3$ of storing a full-rank three-dimensional tensor if $r$ is small compared to $n$.

A compact way of representing $\Abt$ is $\Abt(i_1,i_2,i_3)=\Gc_1(i_1)\Gc_2(i_2)\Gc_3(i_3)$, where the  summation indices are suppressed. Then, for example, for a fixed $i_1, i_2$ application of an $n_3 \times n_3$ SBP operator $D_2$ in the third dimension amounts to multiplying it with the core $\Gc_3(:)$. If we denote the resulting core by $\mathcal{D}_3(i_3)$ we can write the TT approximation of the derivative with respect to the coordinate in the third dimension by 
\[
\widetilde{\bf D}_3(i_1,i_2,i_3)=\Gc_1(i_1)\Gc_2(i_2)\mathcal{D}_3(i_3).
\]
Ignoring boundary conditions and suppressing the spatial dependence, the approximation of the basic update formula in three dimensions with $\Delta u(t)$ denoting the continuous Laplacian
\[
u(t+\Delta t) = 2u(t) - u(t-\Delta t) + (\Delta t)^2\,  \Delta u(t),
 \] 
needed for the wave equation solver, would then be 

\[
\Abt(t+\Delta t) = 2\Abt(t) - \Abt(t-\Delta t)+ (\Delta t)^2 \left[\widetilde{\bf D}_1(t) +\widetilde{\bf D}_2(t) + \widetilde{\bf D}_3(t) \right].
\] 

Just as in the two dimensional case, addition and rounding can be implemented without ever forming the full tensor. We will not describe the details of this procedure and refer to \cite{Oseledets2011} but note that the implementation of the rounding  $\Tc_\delta(\Abt)$ satisfies a relative tolerance criterion
\begin{equation}
\|\Abt-\Tc_\delta(\Abt)\|\leq\frac{\delta}{\sqrt{2}}\|\Abt \|,
\label{eq:TTrounding}
\end{equation}
 where $\sqrt{2}$ comes from dimensional normalization. The cost of performing the addition and rounding of two TT-tensors is $\mathcal{O}(3nr^3)$ but can be reduced to $\mathcal{O}(3nr^2+3r^4)$ by using an intermediate Tucker-expansion.

In this work we use the \texttt{TT-toolbox}\footnote{Available on GitHub at \url{https://github.com/oseledets/TT-Toolbox}.} and the subroutines therein developed by Oseledets et al. This toolbox automatically constructs the TT object class, handles the compression and rounding of TT structures.

\section{The low-rank WaveHoltz method}
\label{sec:LRWHsection}
With the low-rank wave solver in place, we can now generalize the WaveHoltz method outlined in Section~\ref{sec:waveholtz} to the low-rank setting. Again, consider a discretization of $\Omega$ into $m\times m$ blocks $\Omega_{ij}$, we then want to perform blockwise iteration to obtain the WH solution and velocity, denoted by $\Wc_{ij}$ and $\Wc'_{ij}$, respectively. To measure convergence we sum up the norms of the residuals at each block, denoted by $\rho^k_{ij}$, to get the overall measure  $$\rho^k=\sqrt{\sum_{i,j=1}^m(\rho^k_{ij})^2}=\sqrt{\sum_{i,j=1}^m\|\Wc^{k+1}_{ij}-\Wc^k_{ij}\|^2}.$$
Since $\rho^k$ will be bounded by the tolerance used for truncation, we should choose the tolerance $\epsilon^k$ at each iteration informed by the residual. A straightforward way is by scheduling, setting $\epsilon^k=\theta\rho^{k-1}$ for a suitably chosen \textit{scheduling parameter} $\theta\in(0,1)$. The scheduling also works as a control mechanism to keep the rank low, as shown for the Laplace equation in \cite{Appelo2025}. We present the method in Algorithm~\ref{algorithm:LRWaveHoltz} and discuss it below.

\begin{algorithm}
\DontPrintSemicolon
\caption{low-rank WaveHoltz solver (LR-WaveHoltz) }
\label{algorithm:LRWaveHoltz}
\KwIn{Initial data  $\{W^0_{ij}\}_{i,j=1}^m$ and velocity $\{{W^0}'_{ij}\}_{i,j=1}^m$ in low-rank or TT form, Helmholtz frequency $\omega$, discrete Laplacian operators $\{L_{ij}\}_{i,j=1}^m$, discretized forcing in low-rank form $\{F_{ij}\}_{i,j=1}^m$, initial blockwise truncation tolerances $\epsilon^0_{ij}$, convergence tolerance $\epsilon^\star$, grid size $h$.}
\KwOut{Approximate solution and velocity $\Wc^\star,\Wc'^\star$ in low-rank form solving the WaveHoltz fix-point iteration within a tolerance $\epsilon^\star$.}
$k=0$\hfill\Comment{Initialize LRWH iterations}
$\rho^0=1$\hfill\Comment{Initialize LRWH residual}
$T=\frac{2\pi}{\omega}$\hfill\Comment{Calculate fundamental period}
$N_t=\frac{T}{\Delta t}$\hfill\Comment{Calculate number of time steps}
$\epst^0_{ij}=\frac{1}{2N_t}\epsilon^0_{ij}$\hfill\Comment{Initialize tolerance for the wave solver}
\While{$\rho^k>\epsilon^\star$}{
  \For{$ij\in\text{ all multiblocks}$}{
  
        $W^{-1}_{ij}=\Tc^{sum}_{\epst^0_{ij}}[W^0_{ij}-\Delta t{W^0}'_{ij}+\frac{(\Delta t)^2}{2}\Tc^{sum}_{\epst^0_{ij}}(L_{ij}(W^0_{ij}))]$\hfill\Comment{Solution at $-\Delta t$}
        $\Wc_{ij}^0=\frac{3\Delta t}{2T}W^0_{ij}$\hfill\Comment{Initialize WaveHoltz solution at iteration $k$}
        $\Wc_{ij}'^0=0$\hfill\Comment{Initialize WaveHoltz velocity at iteration $k$}
   
    }
    \For{$l=1,\ldots, N_t+1$}{
    $t_l=(l-1)\Delta t$\hfill\Comment{Time stepping}
    \For{$ij\in\text{ all multiblocks}$}{
      
         $[\{W^{l+1}_{ij}\},\{W'^{l+1}_{ij}\}]=\texttt{LRWaveSolver}(\{W^l_{ij}\},\{W'^l_{ij}\},\{F_{ij}\},\Delta t,\epst^k_{ij},t_l)$\\
        $\Wc^{k+1}_{ij}=\Tc_{\epsilon^k_{ij}}^{sum}\bigg(\Wc^{k+1}_{ij}+\frac{2\Delta t}{T}\eta_l(\cos(\omega(t_l+\Delta t))-\frac{1}{4})W^{l+1}_{ij}\bigg)$\\
        $\Wc'^{k+1}_{ij}=\Tc_{\epsilon^k_{ij}}^{sum}\bigg(\Wc'^{k+1}_{ij}+\frac{2\Delta t}{T}\frac{\eta'_l}{2\Delta t}(\cos(\omega t_l)-\frac{1}{4})(W^{l+1}_{ij}-W^{l-1}_{ij})\bigg)$  \\
        }
   
    }
    \For{$ij\in\text{ all multiblocks}$}{
   
         $W^0_{ij}=\Wc^{k+1}_{ij}$\Comment*[r]{Update initial condition for $W_{ij}$}
         $W'^0_{ij}=\Wc'^{k+1}_{ij}$\Comment*[r]{Update initial condition for $W'_{ij}$}
         $\rho^{k+1}_{ij}=\|\Wc^{k+1}_{ij}-\Wc^k_{ij}\|^2$\Comment*[r]{Calculate residual in low-rank form}
        $(\rho^{k+1})^2=(\rho^{k+1})^2+(\rho^{k+1}_{ij})^2$\\
        $\epsilon^{k+1}_{ij}=\max\{K,\theta h\rho^{k+1}_{ij}\}$\\
        $\epst^{k+1}_{ij}=\frac{1}{2N_t}\epsilon^{k+1}_{ij}$
        
    }
    $k=k+1$\hfill\Comment{Update WaveHoltz iteration}
}
\end{algorithm}
 First, we set $\rho^0=1$, the WH period $T$ is obtained from the frequency $\omega$ and the number of time steps $N_t$ are calculated. Moreover, the wave solver tolerance is set $\epst^0=\frac{1}{2N_t}\epsilon^0$ as a scaling of the initial WH tolerance $\epsilon^0$. Then the iteration starts, with each iteration consisting of the same steps. First we initialize the WaveHoltz data using the previous iterates, and calculate the solution at $-\Delta t$. Note that we include the scaling from the WaveHoltz operator in the term $\frac{3\Delta t}{2T}\Wc^0_{ij}$. The current initial data is propagated over the time interval $[0,T]$, adding up the contribution to the integral $\int_0^T(\cos(\omega t)-\frac{1}{4})w(\xb,t)dt$, which is approximated using the trapezoidal rule. The quadrature sum is then truncated using the specified tolerance $\epsilon$. After the time stepping is complete, we store $\Wc^{k+1}$ and $\Wc^k$ to be used for calculating initial data of the next iteration. Then, $\rho^{k+1}_{ij}$ is calculated as $\rho^{k+1}_{ij}=\|\Wc^{k+1}_{ij}-\Wc^k_{ij}\|$, adding it to the total norm of residuals $\rho^{k+1}$. Finally, we then set the block-wise truncation tolerances for the next iteration using scheduling, $\epsilon^{k+1}_{ij}=\theta\rho^{k+1}_{ij}$ and the wave solver tolerance as $\epst^{k+1}_{ij}=\frac{1}{2N_t}\epsilon^{k+1}_{ij}$, where $\theta$ is a specified scheduling parameter.  We then iterate until $\rho^{k+1}$ is smaller than $\epsilon^\star$. 
We note that $\rho^{k+1}_{ij}$ may be very small in blocks before the signal enters. Therefore, truncating with a tolerance proportional to $\rho^{k+1}_{ij}$ may result in an unnecessarily large rank. We therefore propose the alternate strategy of scheduling as $\epsilon^{k+1}_{ij}=\operatorname{max}(K,\theta h\rho^k_{ij})$ for some constant $K>0$. Note that we also include a scaling with the grid size, which is necessary when we want to consider convergence tests.

It now remains to accelerate the two-dimensional LR-WaveHoltz algorithm. For this we use the low-rank Anderson Acelleration method (LRAA) \cite{Appelo2025}. Note that the method outlined in Algorithm~\ref{algorithm:AA} can be extended in a straightforward manner to treat the solution in low-rank form using the truncation operators $\Tc_\epsilon$ and $\Tc_\epsilon^{sum}$. However, since the coefficients $\gamb_k$ should be the same in each block of the discretization and the data now is in matrix form, we need to generalize Problem~\ref{Lemma:AAweights}. The generalization is shown below in Problem~\ref{Lemma:LRAAweights} and the proof can be found in \revone{\ref{proof:LRAA}}. 

\begin{problem}
\label{Lemma:LRAAweights}
 Let $p$ denote the number of blocks in the partitioning of the computational domain, $\|\cdot\|$ the Frobenius norm, $\{\Delta F^{k-j}_i\}_{j=1}^m$ with $\Delta F^{l}_i\in\Rbb^{n\times n}$ the difference matrix in the Anderson acceleration on low-rank form and the set $\{F_i^l\}_{l=1}^m$, $F^k_i\in\Rbb^{n\times n}$ be given. Then, the vector $\gamb^{(k)}$ solving
    $$\gamb^{(k)}=\operatorname{argmin}_{\ub\in\Rbb^m}\bigg(\sum_{l=1}^p\|D^k_l\Lambda(\ub)-F^k_l\|^2\bigg)^\frac{1}{2},$$
    where $\Lambda(\ub)=(\ub^T\otimes I_n)\in\Rbb^{mn\times n}$ and $D^k_l=[\Delta F^{k-1}_l,\ldots,\Delta F^{k-m}_l]\in\Rbb^{n\times mn}$ can be obtained by solving the system $A\gamb^{(k)}=\bb,$ where 
    $$A_{ij}=\sum_{l=1}^p\langle \Delta F^{k-i}_l,\Delta F^{k-j}_l\rangle\quad b_i=\sum_{l=1}^p\langle F^k_l,\Delta F^{k-i}_l\rangle.$$
\end{problem}

We see that the transition from single to multiblock follows as a straightforward sum of the single block terms without changing the overall method. Moreover, the inner products $\langle X,Y\rangle$ can be evaluated efficiently, using only the low-rank factors of $X$ and $Y$ rather than forming an inner product of the full matrices. Finally, we note that $A\gamb^{(k)}=\bb$ is the normal equation for $D^k_l\gamb^{(k)}=F^k_l$. It is well known from standard least squares problems that normal equations often become ill conditioned since the condition number $\kappa$ becomes $\kappa(D^k_l)^2$. We do not see any  problems caused by this in the numerical experiments in Section~\ref{subsec:2Dproblems}.

\section{Numerical Examples}
\label{sec:numericalexamples}

Unless specified otherwise, we discretize in space using the 4th order SBP-FD operators by Mattsson in \cite{Mattsson2004} and set the penalty parameter to $\tau=15$. For the time discretization, we use the leap-frog method as mentioned in the previous sections with the time step $\Delta t=0.15 h$.The residuals in LRWH are calculated using the Frobenius norm, while the errors will be made grid-independent by scaling the Frobenius norm with $h$.Recall that the  truncated SVD operator $\Tc_\epsilon$  approximates a matrix $W$ within a given accuracy $\epsilon$ as $\|W-\Tc_\epsilon(W)\|<\epsilon$. Thus to compute the error in the continuous 2-norm we have
\[
\left( \int \int (e(x,y))^2 dx dy \right)^{1/2} \approx  h \| E \|.   
\]
In other words we should take $\epsilon = \frac{{\rm TOL}}{h}$ in order to get an error ${\rm TOL}$ in the continuous 2-norm.

\begin{figure}[]
    \centering
    
    \begin{minipage}[b]{0.52\textwidth}
        \centering
        \includegraphics[width=\textwidth]{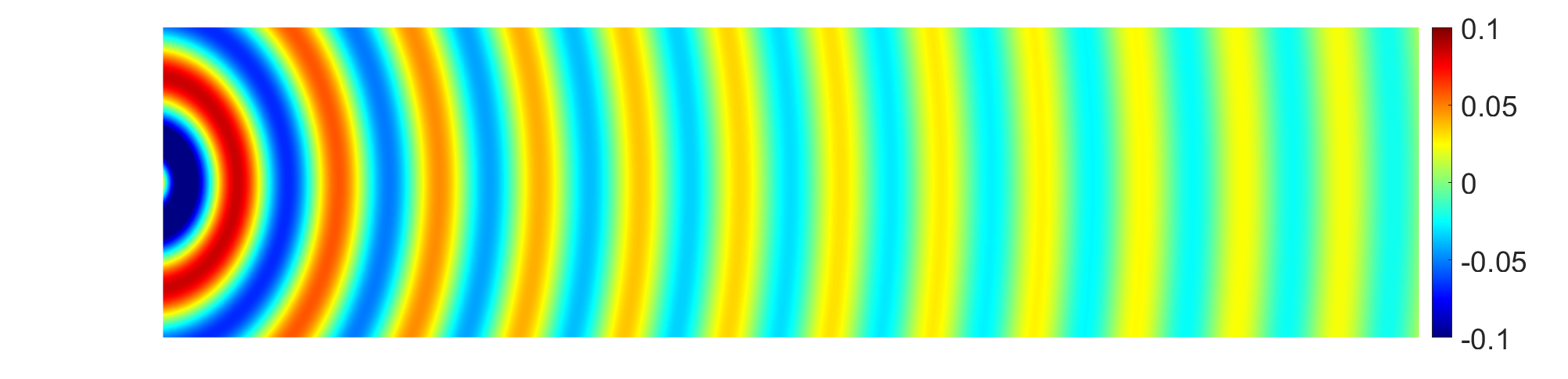} %
        \captionof{figure}{The free-space Green's function $G(x,y)$ of the Helmholtz equation centered at $(-0.1,0.5)$ throughout the domain $[0,5]\times[0,1]$.}
        \label{fig:ex1GreensFunction}
    \end{minipage}
     \hspace{0.01\textwidth}
    \vrule width 1pt
     \hspace{0.01\textwidth}
    \begin{minipage}[b]{0.4\textwidth}
        \centering
        \begin{tabular}{|c|c|}
            \hline
             Degrees of freedom & PPW \\
            \hline
            $51^2$ & $20$ \\
            $101^2$ & $40$ \\
            $201^2$ & $80$ \\
            $401^2$ & $160$ \\
            $801^2$ & $320$\\
            $1601^2$ & $640$\\
            \hline
        \end{tabular}
        \captionof{table}{Number of degrees of freedom per block and the resulting points per wavelength using a frequency $\omega=5\pi$ and unit wave speed.}
        \label{table:Ex1PPW}
    \end{minipage}
\end{figure}
\subsection{Problems in 2D}
\label{subsec:2Dproblems}

\subsubsection{Timing and compression tests} 
We begin by demonstrating the efficiency of the low-rank framework by comparing runtimes  and compression for the Green's function. Let $\Omega=[0,5]\times[0,1]$. We set the wave speed to $c=1$ and let the forcing be given by a point source located at $(x_0,y_0)$ oscillating with a frequency $\omega$.  \revtwo{The solution is given by the Green's function 
\begin{equation}
    G(x,y)=\frac{i}{4}H_0^{(1)}(\omega r(x,y)),\quad r(x,y)=\sqrt{(x-x_0)^2+(y-y_0)^2},
    \label{eq:Greensfnc}
\end{equation}
where $r(x,y)$ denotes the distance from the evaluation point $(x,y)$ to the point source and $H_0^{(1)}$ denotes the Hankel function of the first kind.} 

To quantify the cost of low-rank operations, we time the SBP approximation to the Laplacian of the Green's function followed by an addition and application of the truncation operator $\Tc_\epsilon$ to mimic a time step. \revtwo{To get stable timings we repeat this 100 times. We fix the frequency to be $\omega=5\pi$ and the location of the point source to be $(x_0,y_0)=(-0.05,0.5)$. The domain $\Omega$ is partitioned into $5 \times 1$ blocks and vary the number of points per wavelength (PPW) (see Table~\ref{table:Ex1PPW}). To compare the full- and low-rank methods, we consider truncation tolerances: $\epsilon=10^{-2}h^{-1}$, $\epsilon=10^{-3}h^{-1}$ and $\epsilon=10^{-4}h^{-1}$.} 

\begin{figure}[]
        \centering
        \includegraphics[width=0.99\linewidth]{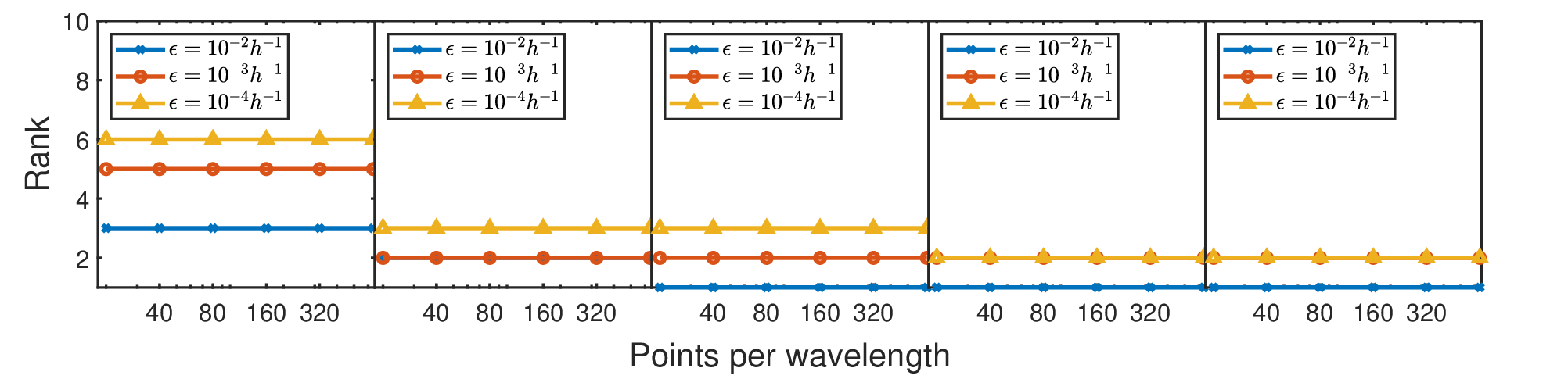}    
       \caption{Ranks of the truncated Greens function $\Tc_\epsilon(G(X,Y))$ in each of the $5\times 1$ blocks using the truncation tolerances \revtwo{$\epsilon=10^{-2}h^{-1},\epsilon=10^{-3}h^{-1}$ and $\epsilon=10^{-4}h^{-1}$} for varying number of points per wavelength.}
      \label{fig:Greenranks}
\end{figure}

We investigate how the rank of the Green's function behaves for various levels of truncation. Since the rank is an inherent property of the function itself, this will work as a guideline on how to choose the tolerance depending on the desired rank.

\begin{figure}[]
    \centering
    \includegraphics[width=0.99\linewidth]{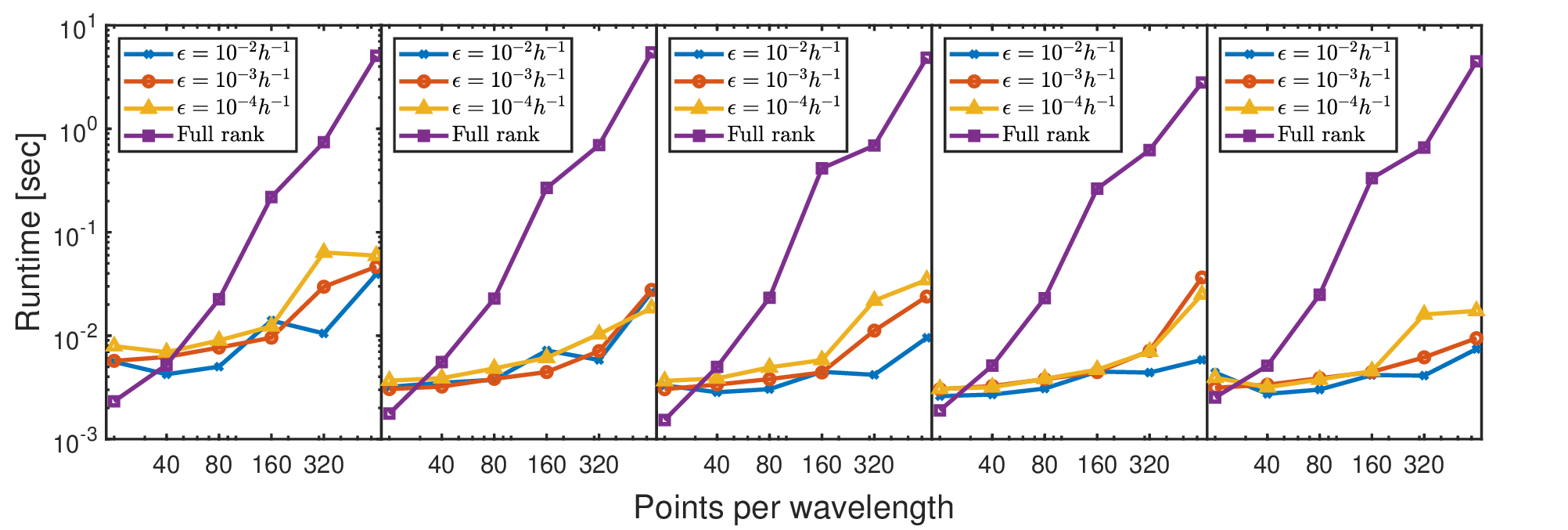}
    \includegraphics[width=0.99\linewidth]{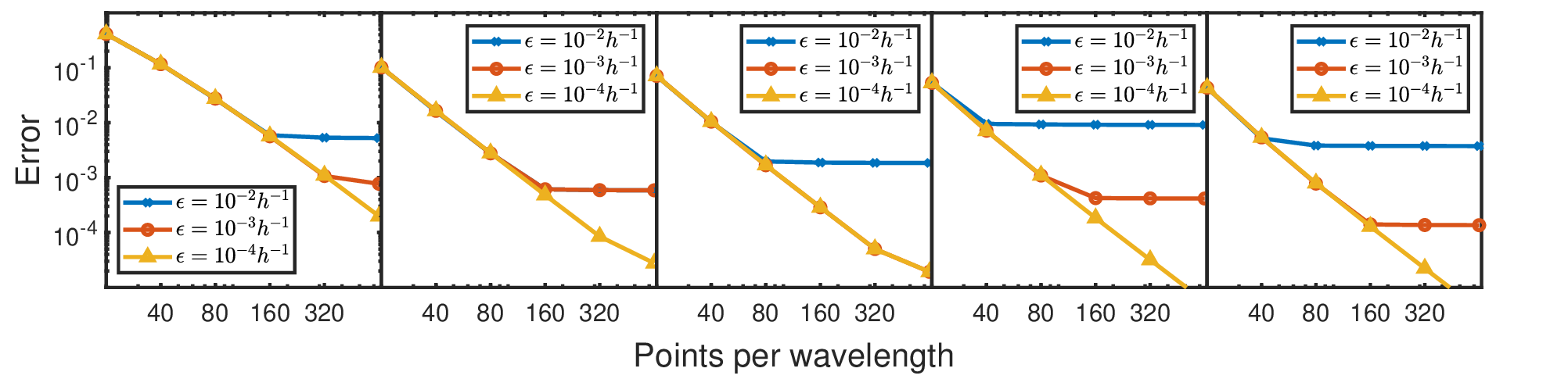}
    \caption{Time to simulate hundred time steps in each block using full- and low-rank methods (top) \revtwo{discrete errors between analytical and numerical Laplacian as a function of points per wavelength using the low-rank discretization with truncation tolerances $\epsilon=10^{-2}h^{-1},\epsilon=10^{-3}h^{-1}$ and $\epsilon=10^{-4}h^{-1}$.}}
    \label{fig:GreensRuntime}
\end{figure}

We illustrate the Greens function in Figure~\ref{fig:ex1GreensFunction} and the ranks of $\Tc_\epsilon(G(X,Y))$ in blockwise order in Figure~\ref{fig:Greenranks} for the different choices of truncation tolerances. We note that the rank is higher in the leftmost blocks, which is expected since these are close to the source. The rank then decreases as we move away from the point source as expected.

The run-times, ordered blockwise, are shown in Figure~\ref{fig:GreensRuntime} together with the approximation error for the different truncation tolerances. The low-rank method is more efficient for larger numbers of PPW and further from the source, with a potential gain of more than one order of magnitude but remains efficient even for smaller PPW. Notably, the difference in runtime is small when decreasing $\epsilon$. The runtime also decreases as we move away from the point source, agreeing with the results in Figure~\ref{fig:Greenranks}. The errors stagnate at higher level for higher tolerances in each block. We highlight that we can gain orders of magnitude in error while maintaining a small runtime even for the smallest tolerance. It is also clear that the overall magnitude of errors decrease further away from the source as the norm of $G(X,Y)$ decreases.

\subsubsection{Acceleration of free-space and two-corner problems}
We now turn to the full LR-WaveHoltz (LRWH) method in two dimensions and consider the domain $\Omega=[0,2]\times [0,1]$ using the method in Algorithm~\ref{algorithm:LRWaveHoltz}. First, we investigate if how well the low-rank Anderson Acceleration method with memory parameter $M$ ($\LRAAM$) works for two problems. Since the WaveHoltz method is known to converge slower for problems with (partially) trapped waves \cite{ROTEM2026} we compare problem with homogeneous Neumann boundary conditions on three sides, with a free-space problem. We partition the domain into $6\times 3$ blocks with a point source located at $(x_0,y_0)=(0.1,0.5)$ approximated by 
\begin{equation}
f(x,y)=-\frac{1}{\delta^2}\operatorname{exp}\bigg(-\frac{(x-x_0)^2+(y-y_0)^2}{\delta^2}\bigg),\quad\delta=\frac{1}{2\omega},
\label{eq:Gaussian}
\end{equation}
and set the frequency to $5\pi$. Let $X(i,j)=ih$ and $Y(i,j)=jh$ denote the grid points in a block. Then $f(X,Y)=U_fS_fV_f^T$, where 
\[
U_f=\operatorname{exp}(-(X(:,1)-x_0)^2/\delta^2)\quad S_f=-1/\delta^2\quad V_f=\operatorname{exp}(-(Y(1,:)-y_0)^2/\delta^2).
\]

Finally, we also have to set the stopping tolerance $\epsilon^\star$. We discretize using $n=26$ gridpoints per direction in each block so that  the \revone{mesh size} is $h=\frac{1}{75}$ and $PPW=30$. Since we use fourth order operators we use the rule-of-thumb presented in \cite{Appelo2025Rule} and choose $\epsilon^\star=10^{-3}$ so that the error in integrated norm becomes $h\epsilon^\star=1.33\times 10^{-5}$. 

\begin{figure}[]
    \centering
    \includegraphics[trim={2cm 0 2cm 0}, clip, width=\linewidth]{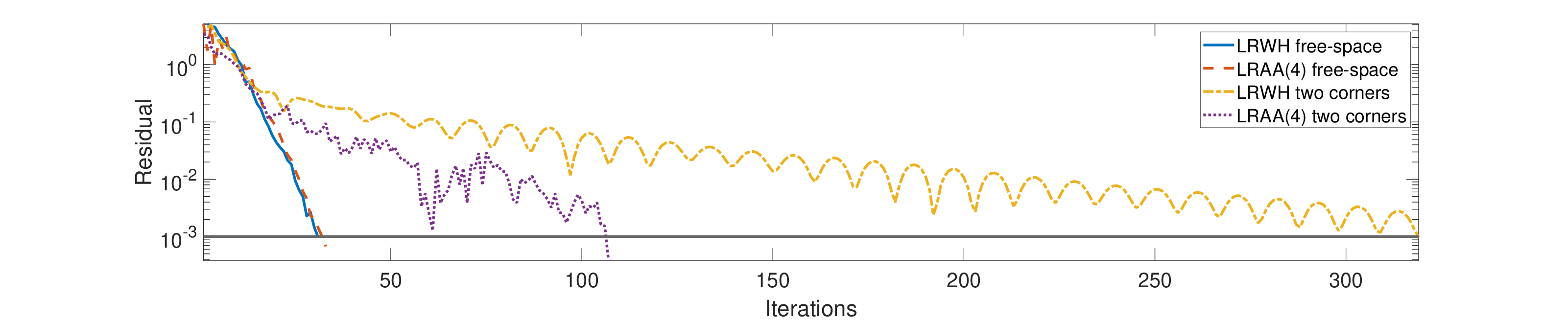}
    \caption{Residuals obtained using the LRWH and $\operatorname{LRAA}(4)$ methods for a free-space problem and problem with two corners, respectively. This is obtained using the stopping tolerance $\epsilon^\star=10^{-3}$ and scheduling.}
    \label{fig:AccelerationEx}
\end{figure}

We solve using both the LRWH and $\operatorname{LRAA}(4)$ method. Note that we use scheduling in $\LRAAM$ method, described in more  detail in Section~\ref{subsubsec:freespaceproblem}. The resulting residuals are displayed in Figure~\ref{fig:AccelerationEx}. The convergence for the problem with corners is slower than that of the free-space problem, as expected and the LRAA method outperforms LRWH for the problem with corners. 
\subsubsection{Free-space problem}
\label{subsubsec:freespaceproblem}
\begin{figure}[]
    \centering
    \includegraphics[width=0.5\linewidth]{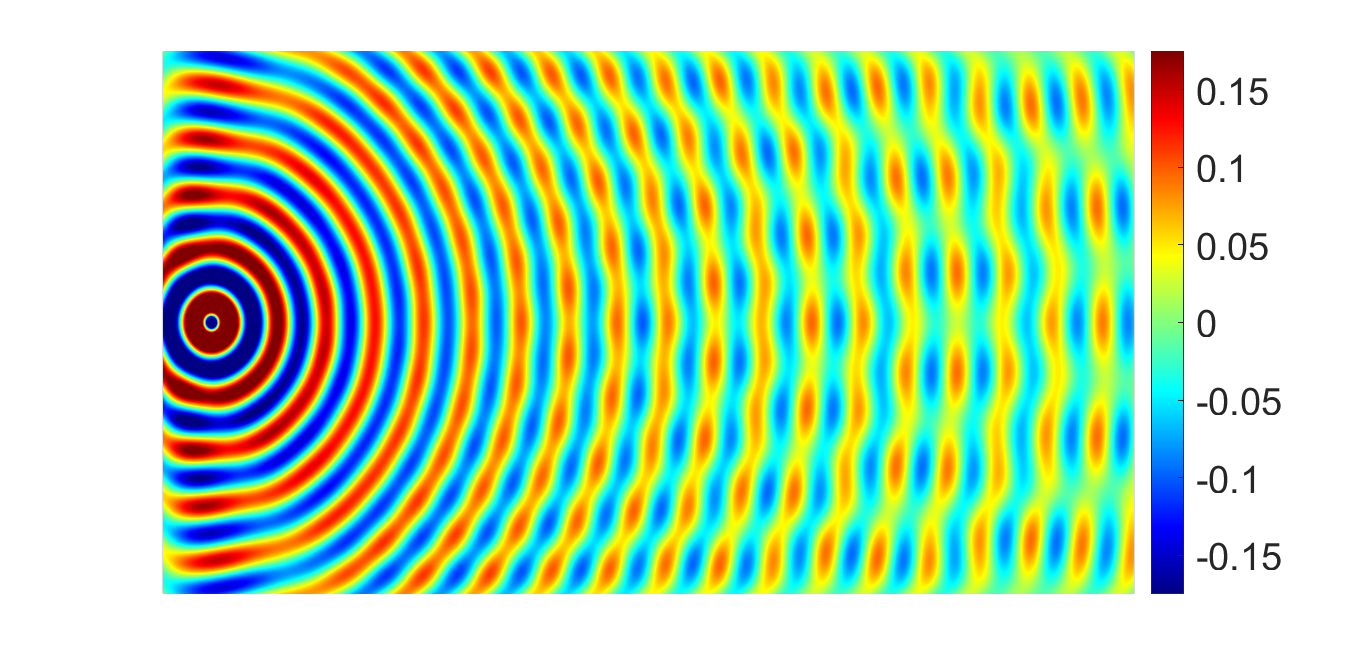}
    \includegraphics[width=0.45\linewidth]{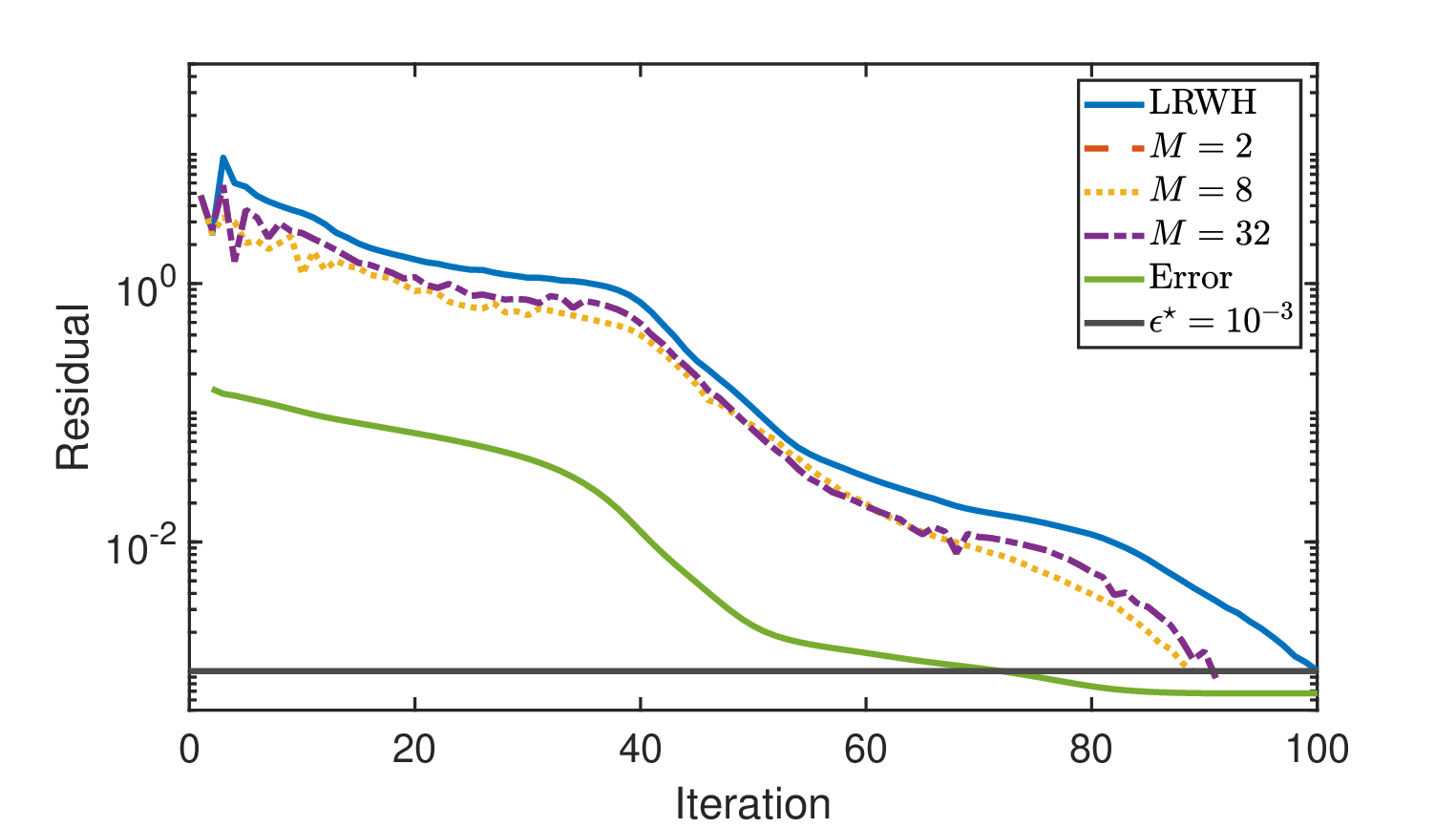}
    \caption{Free-space problem. The resulting Helmholtz solution in a domain partitioned into $6\times 3$ blocks caused by a point source when imposing outflow conditions along all boundaries (left) residuals for LRWH and $\min(\rho_G,\rho_X)$ obtained from LRAA using memory parameters $M=2,8,32$ with stopping tolerance $\epsilon^\star=10^{-3}$ and scheduling parameter $\theta=1$ together with the discrete error in (right).}
    \label{fig:OutflowEx}
\end{figure}

Motivated by the results in the previous section, we carry out a more detailed study of the free-space problem. For the spatial discretization we instead use $n=101$ grid points per direction in each block with frequency $\omega=20\pi$, maintaining $PPW=30$ with the new mesh size $h=\frac{1}{300}$. We consider both the LRWH and $\LRAAM$ method. In the accelerated case, we need to make a decision as to what residual should be used to measure convergence. First, there is the error between the iterate and fixed point operator $\rho_G=\|X^k-G(X^k)\|$, and the error between consecutive iterates $\rho_X=\|X^k-X^{k+1}\|$. Since the truncation should be with respect to the smallest residual we use a strategy to choose $\epsilon=\max\{K,\theta h\min\{\rho_X,\rho_Y\}\}$, where $\theta>0$ is the scheduling parameter and $K>0$ a constant. When the computational domain is covered by a single grid the first term is not needed, but when the domain is discretized by multiple block and the solution in some of the blocks is the same size as machine precision the first term will prevent artificial rank inflation early on in the iteration. Once the solution enters the block the first term is effectively removed from the computation of $\epsilon$. 

\begin{figure}[]
\centering
\includegraphics[width=0.9\linewidth]{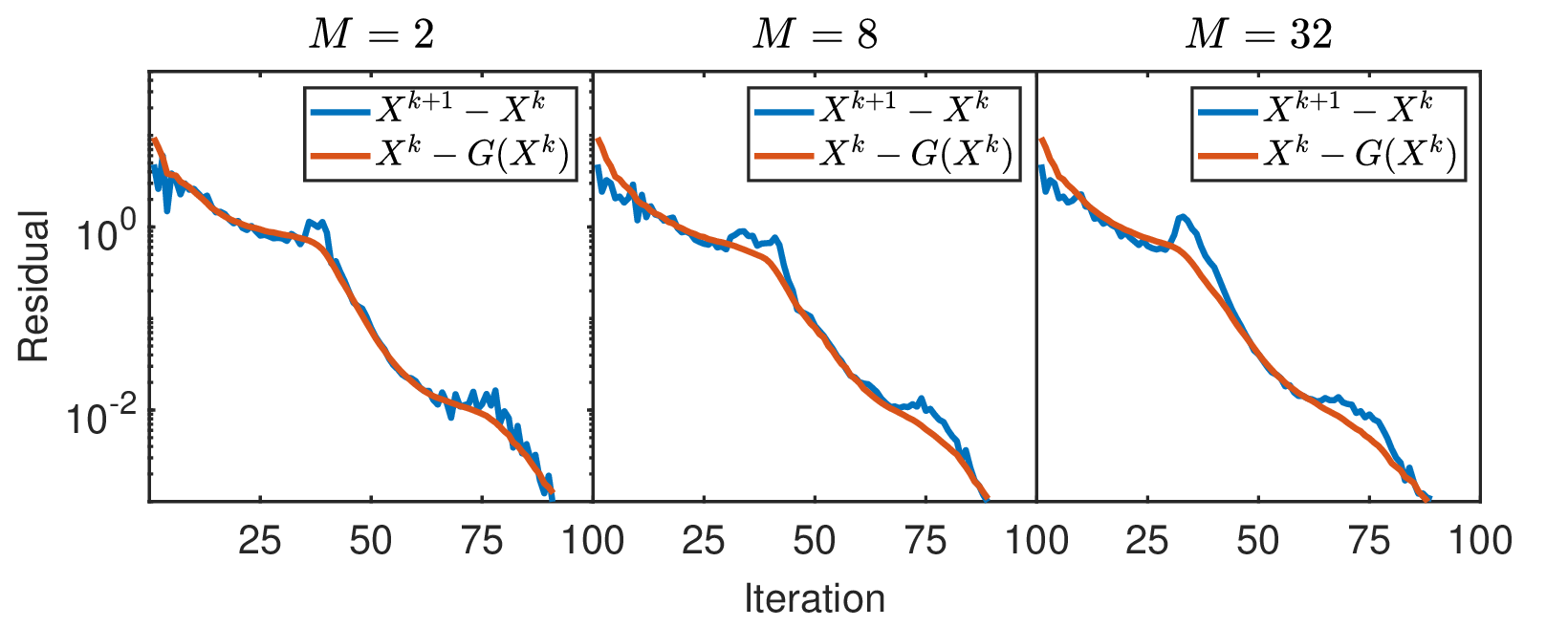}
\caption{Free-space problem. The total AA(M) residuals $\rho_G$ and $\rho_X$ for increasing memory parameters $M=2,8,32$ obtained using the truncation tolerance $\epsilon=\max\{K,\theta h\min\{\rho_X,\rho_Y\}\}$ with scheduling parameter $\theta=1$ and $K=10^{-5}$.}
\label{fig:ChannelRes}
\end{figure}

The resulting solution is presented in Figure~\ref{fig:OutflowEx} together with the LRWH and $\LRAAM$ residuals for $M=2,8,32$ as well as the error between the LRWH method and full rank solution obtained using SBP operators to discretize \eqref{eq:modelproblem} and then performing an explicit inverse. We see a clear convergence, but no significant acceleration using the $\LRAAM$ method. It is also clear that the discrete error is about the size of the residual at convergence. In Figure~\eqref{fig:ChannelRes}, we present the residuals $\rho_X$ and $\rho_G$ for different $M$. The residuals behave similarly for all values of $M$. It is not clear that either of the residuals always dominates the other one, confirming that we should use $\min(\rho_G,\rho_X)$ to adjust the truncation tolerance.  It is, however, clear that $\rho_G$ is smoother than $\rho_X$, similar to what is reported in \cite{Appelo2025}. In all examples, the stopping tolerance is set to $\epsilon^\star=10^{-3}$, corresponding to an integrated residual of $h\epsilon^\star=3.33\times 10^{-6}$ at convergence. 

\begin{figure}[]
\centering
   \includegraphics[width=\linewidth]{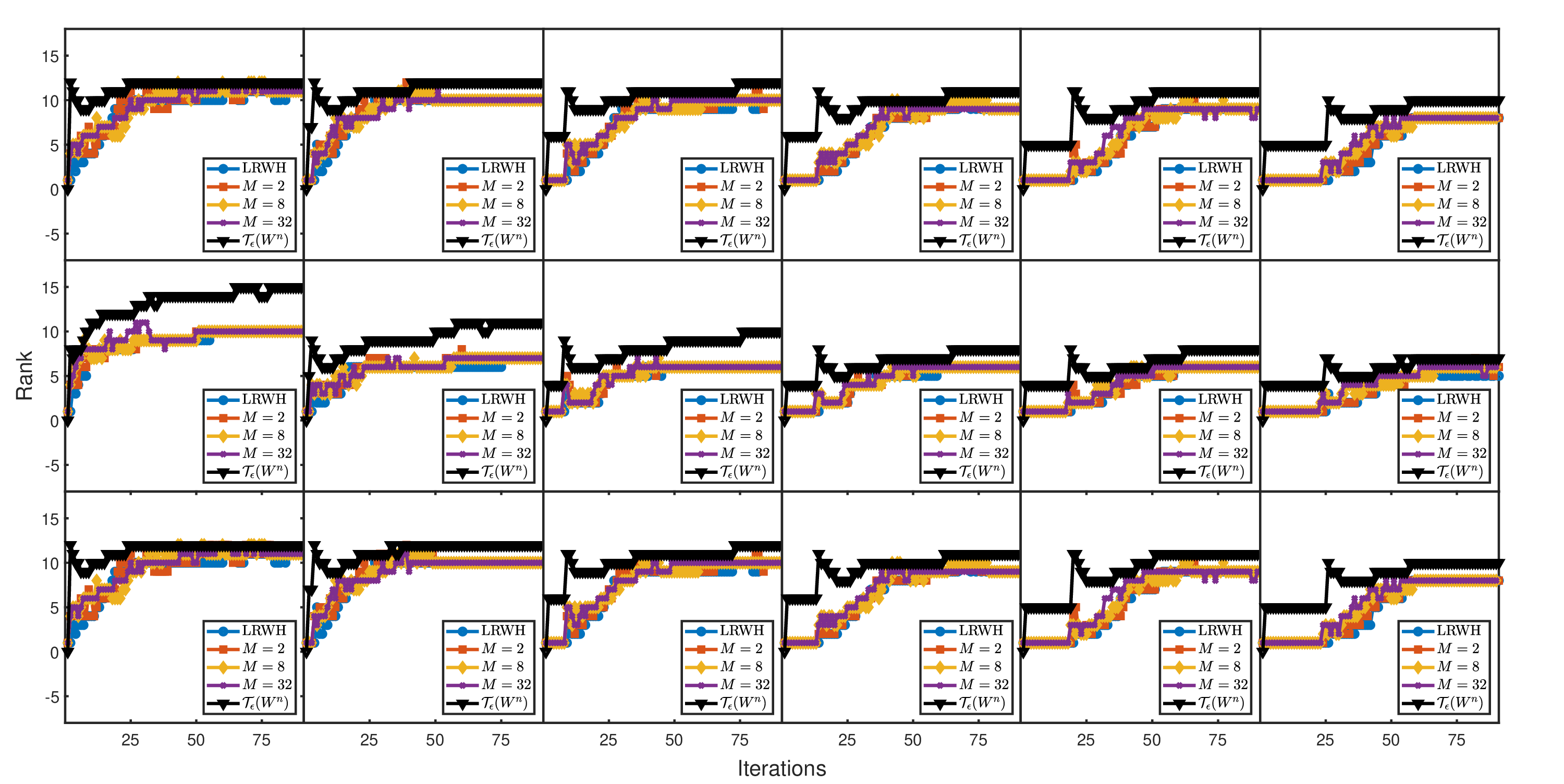}
\caption{Free-space problem. Distribution of ranks for the LRWH and LRAA methods in each of the $6\times 3$ blocks using the scheduling parameter $\theta=1$ and stopping tolerance $\epsilon^\star=10^{-3}$ for the memory parameters $M=2,8,16$ as well as the rank of the full rank solution $W^n$ truncated at the tolerance at the $n$'th iteration.}
\label{fig:ChannelRank}
\end{figure}

In Figure~\ref{fig:ChannelRank}, we present the ranks for both the LRWH and LRAA(M) methods, as well as the rank obtained by truncating the full-rank solution $W$ using the blockwise truncation tolerance $\epsilon^k_{ij}$ for each iteration $k$, denoted by $\Tc_{\epsilon^k}(W)$. In all blocks we see that the numerical rank remains bounded by that of $\Tc_{\epsilon^k}(W)$. Moreover, the ranks for both the LRWH and LRAA(M) methods remain very similar and show a monotonic growth.Just as for the Green's function the rank of the solution decreases with the distance from the source. For this problem $\LRAAM$ (for all choices of memory parameter) and LRWH performs similarly well.
\revonesecondgo{\subsubsection{Convergence of the WaveHoltz method}
\begin{figure}[]
    \centering
    \includegraphics[width=0.4\linewidth]{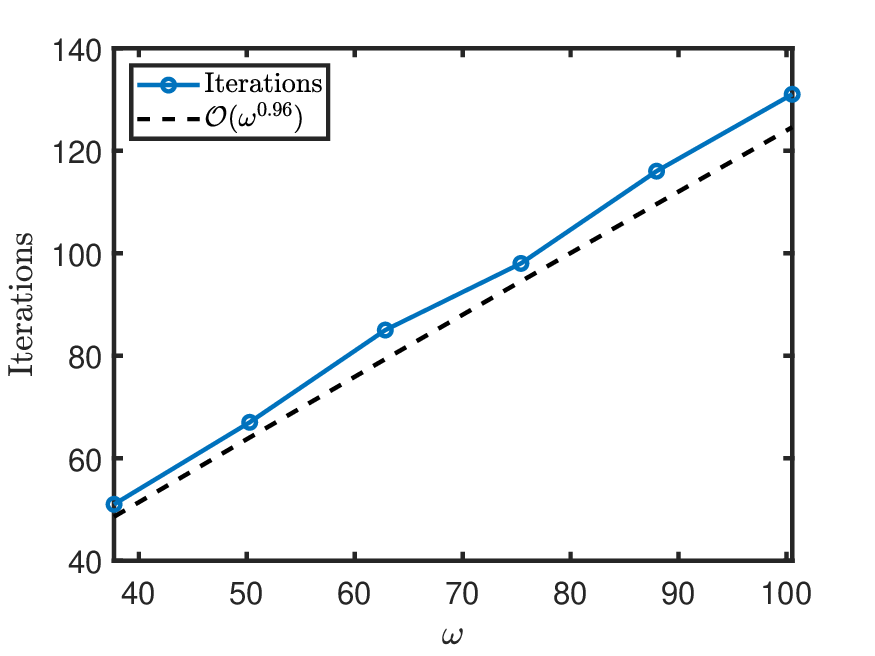}
   \includegraphics[width=0.4\linewidth]{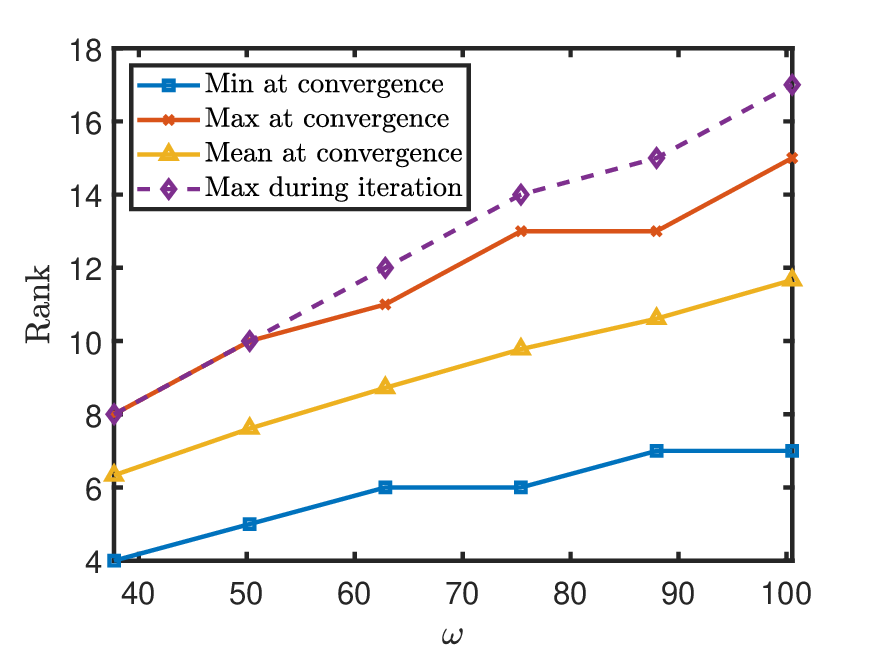}
    \caption{Free-space problem. The number of iterations required for LRAA(2) to converge with stopping tolerance $\epsilon^\star=10^{-3}$ and scheduling parameter $\theta=0.5$ together with the power-law fit (left) and the minimal-, maximal- and mean rank throughout all of the 18 blocks at convergence, together with the maximal rank during the iterations (right).}
    \label{fig:WHConvergence}
\end{figure}

To consider how the acceleration and truncation affects the overall WaveHoltz method we now investigate the evolution of the numerical rank during the iteration and the number of iterations required to converge. We fix $\PPW=10$ and consider the free-space problem above while varying the frequency $\omega$.

The results are presented in Figure~\ref{fig:WHConvergence}. In the left figure we see that the number of iterations scales like $\mathcal{O}(\omega^{0.96})$, which agrees very well with  analysis in \cite{ROTEM2026}. The ranks are presented in the right figure, demonstrating an increase of rank as the frequency grows. However, we see that there is essentially no intermediate rank inflation as the frequency grows.} 
\subsubsection{Half-space problem}
\begin{figure}[]
    \centering
    \includegraphics[width=0.5\linewidth]{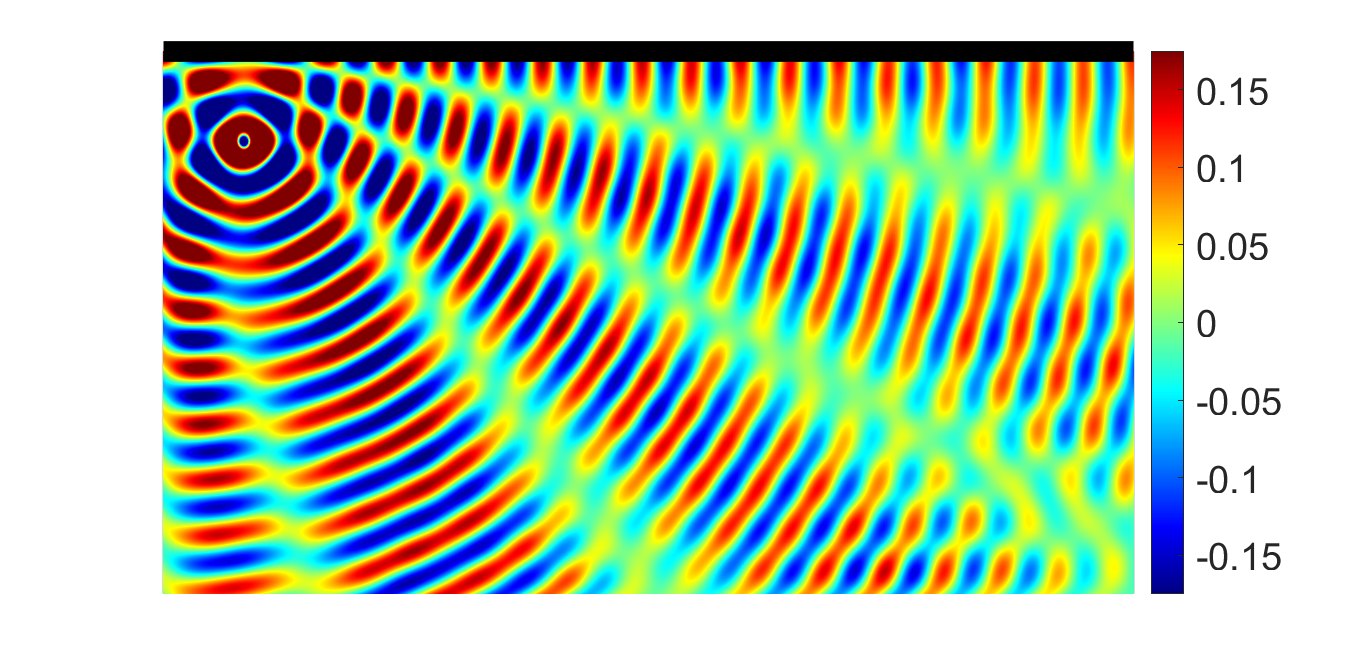}
   \includegraphics[width=0.45\linewidth]{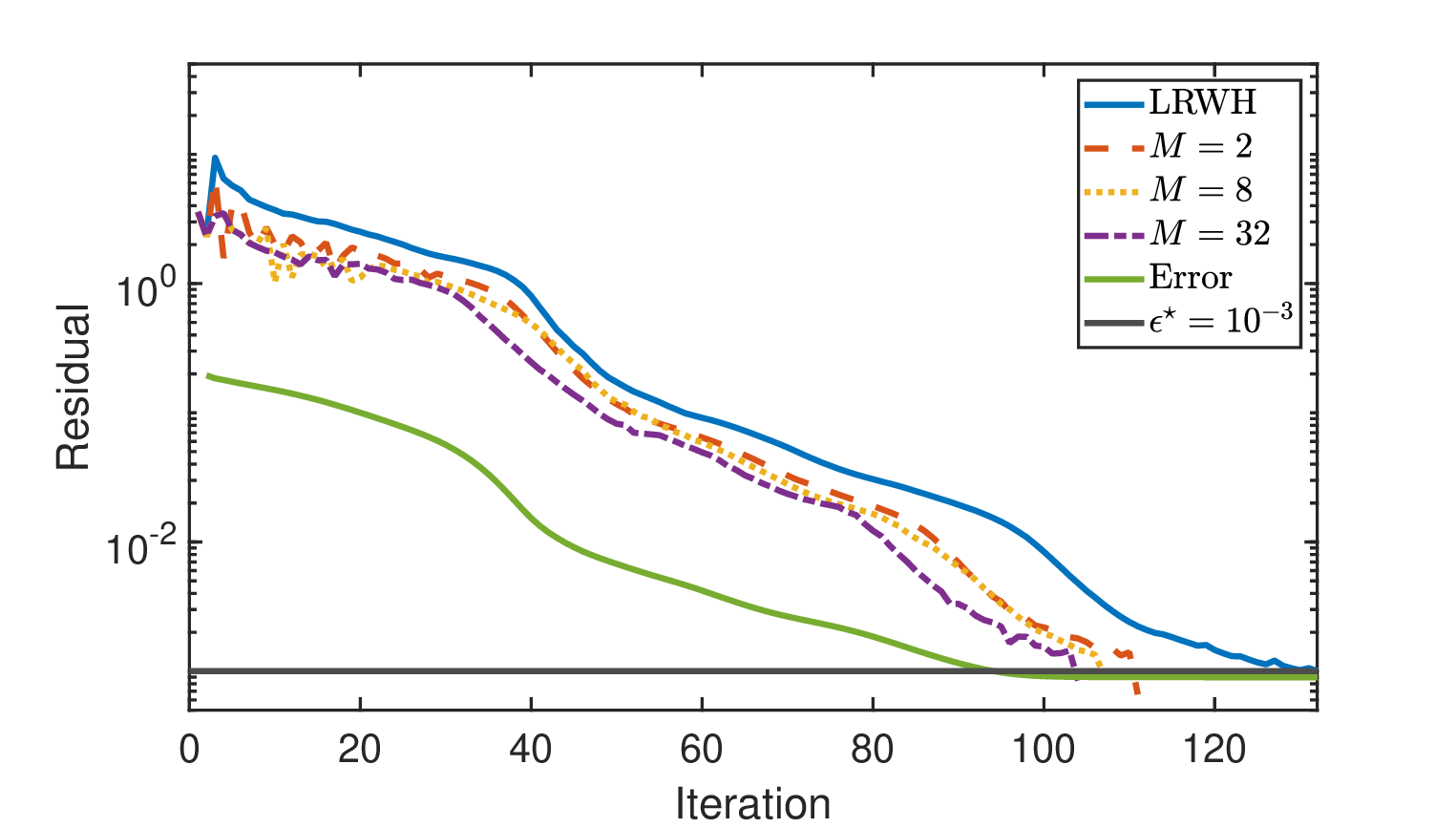}
    \caption{Half-space problem. The resulting Helmholtz solution in a domain partitioned into $3\times 6$ blocks generated by a point source below a reflecting water surface  (left) residuals for the LRWH and LRAA methods using memory parameters $M=2,8,32$ with stopping tolerance $\epsilon^\star=10^{-3}$ and scheduling parameter $\theta=1$ together with the discrete error (right).}
    \label{fig:ReflectingEx}
\end{figure}
We now turn to the situation where the point source is located just below the water surface. Assuming a still water with slow wind speeds above the surface allows us to approximate the surface with a straight line. Since the signal will reflect at the surface, the boundary conditions must be changed to a Neumann condition, indicated by the solid line in Figure~\ref{fig:ReflectingEx}. The signal is still modeled by \eqref{eq:Gaussian}, but moved to $(x_0,y_0)=(1/6,5/6)$ and the domain is partitioned using a total of eighteen blocks, three in the vertical direction and six in the horizontal direction. Each block is then discretized using $101^2$ degrees of freedom so that $\PPW=30$ \revone{and $h=\frac{1}{300}$}. We use the same method as in the free-space example to determine $\epsilon^k$. The resulting solution is shown in Figure~\ref{fig:ReflectingEx} together with the LRWH and $\LRAAM$ residuals for memory parameters $M=2,8,32$, as well as the error calculated in Frobenius norm. The solution behaves as expected, and both the LRWH and $\LRAAM$ methods show a clear convergence. We again note that there only is a slight reduction in the number of iterations between the LRWH and $\LRAAM$ method, but increasing $M$ does not have any significant effect on the number of iterations. The discrete error is about the size of the residual at convergence.

\begin{figure}[]
\centering
   \includegraphics[width=\linewidth]{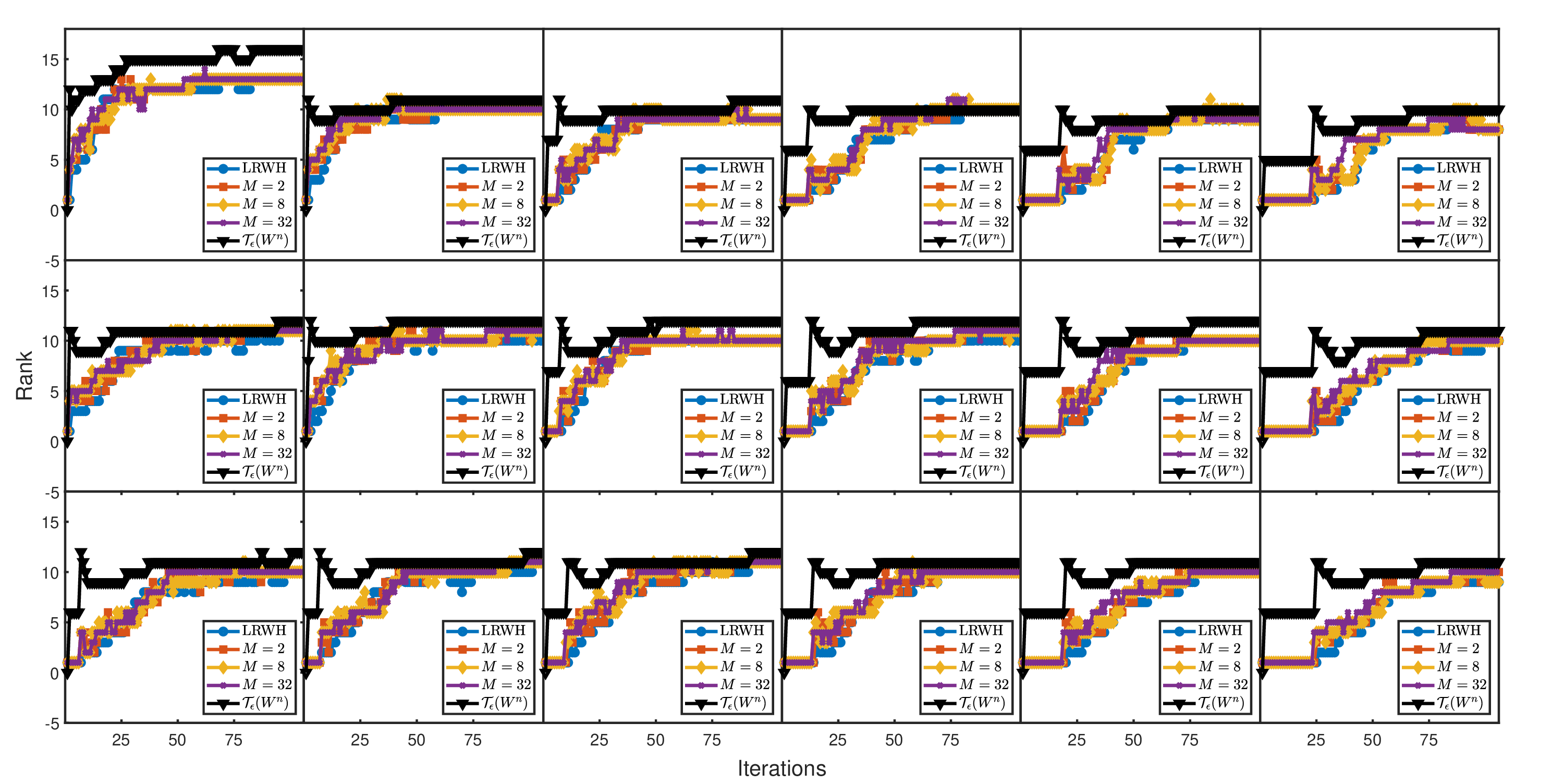}
\caption{Half-space problem. Distribution of ranks for the LRWH and LRAA methods in each of the $6\times 3$ blocks using the scheduling parameter $\theta=1$ and stopping tolerance $\epsilon^\star=10^{-3}$ for the memory parameters $M=2,8,16$ as well as the rank obtained from the full rank solution $W^k$ truncated at the tolerance at the $k$'th iteration.}
\label{fig:ReflectionRanks}
\end{figure}

The resulting ranks are presented in Figure~\ref{fig:ReflectionRanks} together with the rank of $\Tc_{\epsilon^k}(W)$. Similar to the previous example we note that the numerical rank remains bounded by the rank of $\Tc_{\epsilon^k}(W)$ and demonstrates a nearly monotonic growth. The overall size of the final ranks also follow the expected behavior with a large rank close to the pointsource.

\subsubsection{Layered free-space problem}
\begin{figure}[]
    \centering
    \includegraphics[width=0.55\linewidth]{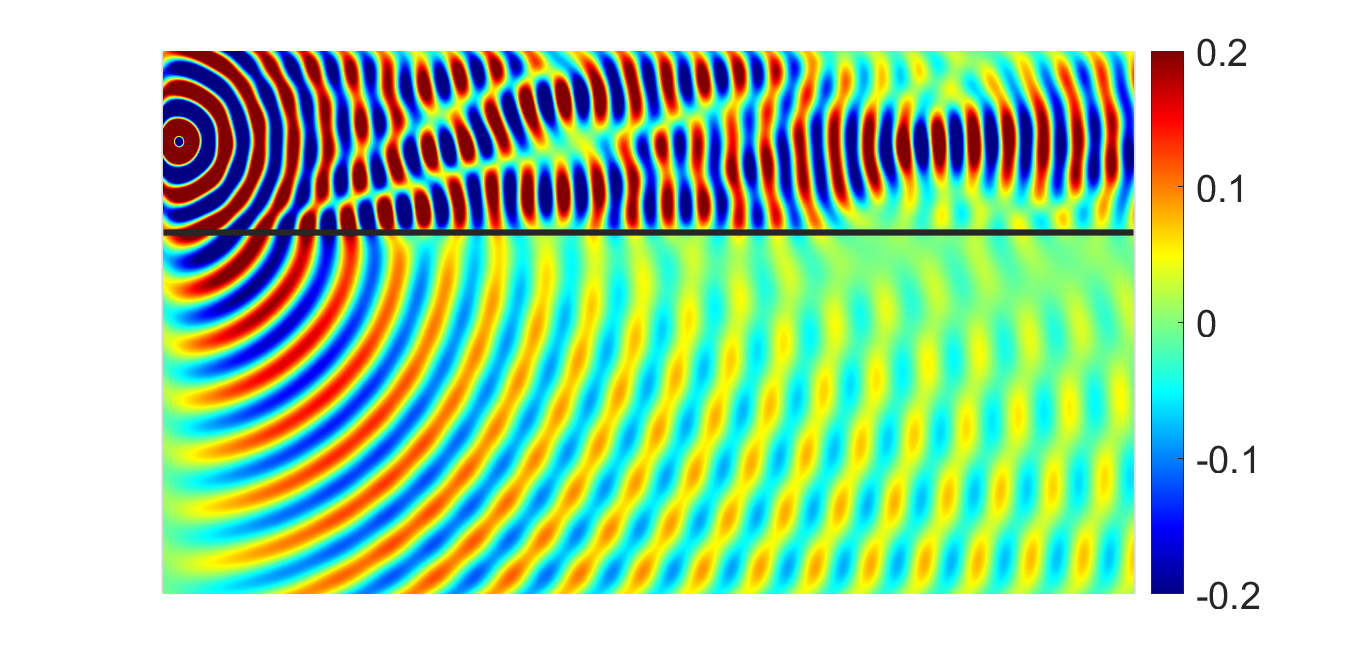}
   \includegraphics[width=0.4\linewidth]{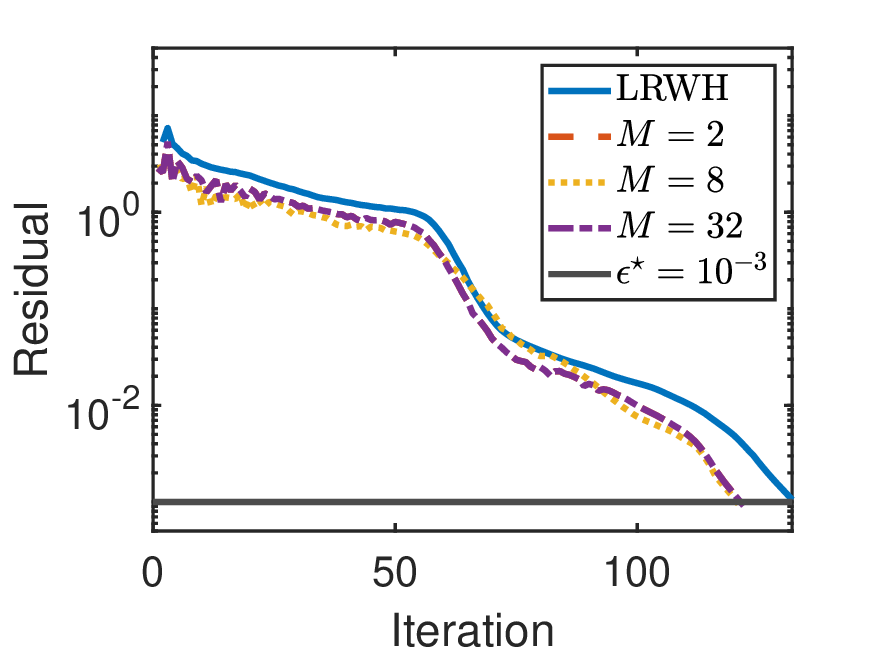}
    \caption{Layered free-space problem. The resulting Helmholtz solution in a discontinuous medium partitioned into $6\times 3$ blocks caused by a point source located above the interface (left) residual for the LRWH and $\operatorname{min}(\rho_G,\rho_X)$ obtained from LRAA using memory parameters $M=2,8,32$ with stopping tolerance $\epsilon^\star=10^{-3}$ and scheduling parameter $\theta=1$ (right).}
    \label{fig:DiscontinuousEx}
\end{figure}

\begin{figure}[]
    \centering
    \includegraphics[width=\linewidth]{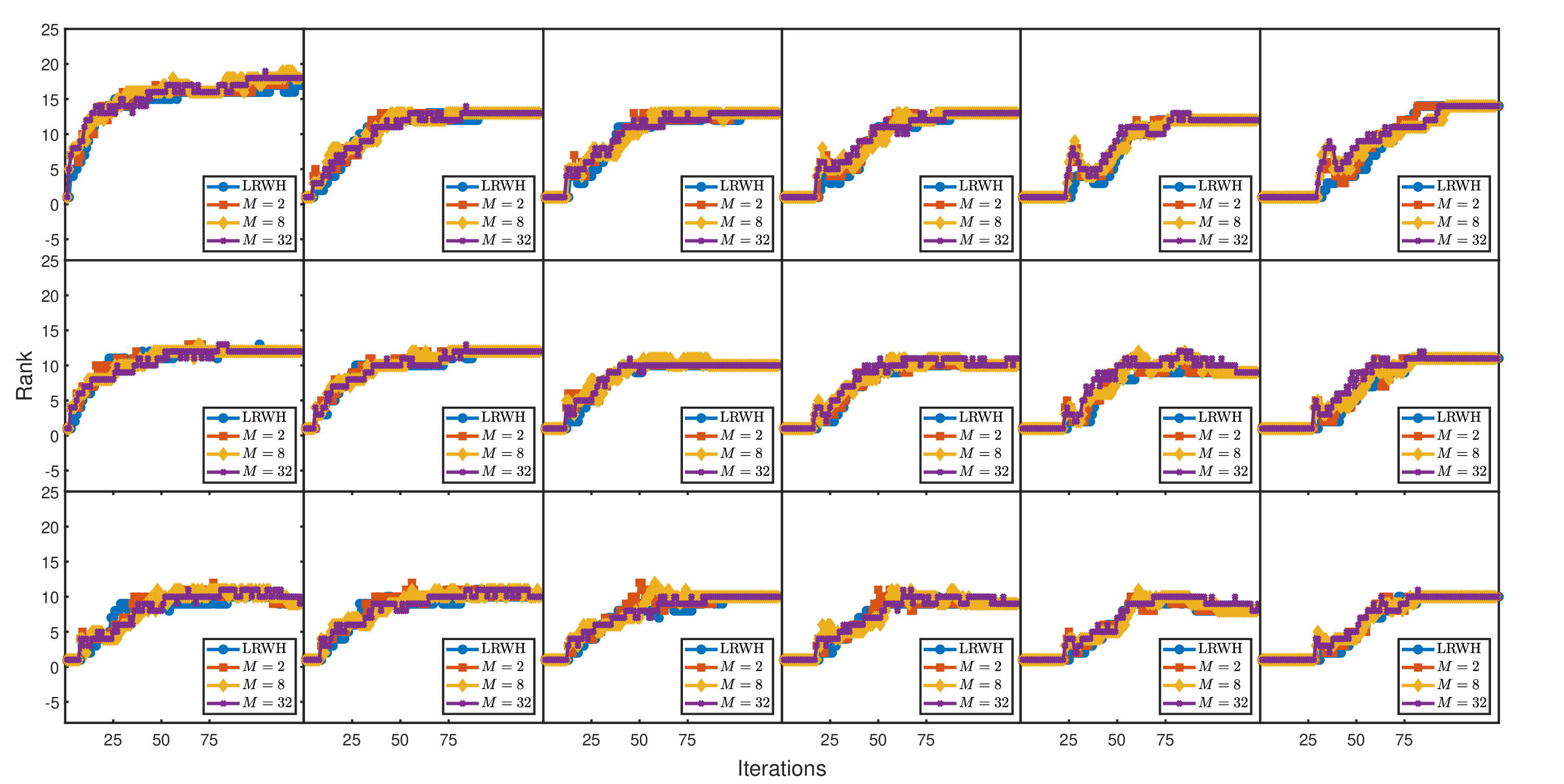}
\caption{Layered free-space problem. Distribution of ranks for the LRWH and LRAA methods in each of the $6\times 3$ blocks using the scheduling parameter $\theta=1$ and stopping tolerance $\epsilon^\star=10^{-3}$ for the memory parameters $M=2,8,16$ as well as the rank obtained from the full rank solution $W^k$ truncated at the tolerance at the $k$'th iteration. }
    \label{fig:DiscontinuousRanks}
\end{figure}

Another type of common problem is when the wave speed varies discontinuously between two neighboring computational domains. This may be because of local layerings in the water or because of interactions with a fluidic seafloor, as is the case in sediments. \change{To model this situation, we partition the domain as $\Omega=\Omega_1\cup\Omega_2$ into a water domain $\Omega_1$, a sediment domain $\Omega_2$ and denote the wave speed in the respective domains as $c_1$ and $c_2$.}\revone{ We discretize each block by $101^2$ degrees of freedom so that $h=\frac{1}{300}$ and $\text{PPW}=30$}. If the point source is located in $\Omega_1$ there will be reflections from the interface between the two domains. If the wave speed is higher in the sea floor than in the water, i.e. $c_2>c_1$, it follows from Snell's law that there exists a critical angle $\theta_c=\operatorname{arccos}(c_1/c_2)$ such that all incoming waves propagating within the aperture $2\theta_c$ does not loose any energy due to sea floor interactions. This mechanism is what explains the emergent waveguide behavior in ocean acoustics, particularly common when analyzing continental shelves, where $\theta_c$ usually lies between $10^\circ$ and $30^\circ$ \cite{Jenssen2011}. We take inspiration from these numbers and set $c_1=0.7, c_2=1$. Using the same method for choosing the truncation tolerance as in the previous examples, we obtain the solution and residuals shown in Figure~\ref{fig:DiscontinuousEx}. The solution clearly demonstrates a layering effect away from the point source, while letting signals through at a narrower angle closer to the source. The residual demonstrates convergence both of the LRWH and $\LRAAM$ methods for $M=2,8,32$ but see no clear acceleration when increasing $M$. In Figure~\ref{fig:DiscontinuousRanks} we also see the expected growth in the ranks, being nearly monotonic for LRWH and all choices of memory parameter. Note the small peaks in the ranks of the rightmost blocks. This is an example of the solution having a small residual in the block, rendering the truncation tolerance unnecessarily small.  

\subsubsection{Layered half-space problem}
\begin{figure}[]
    \centering
    \includegraphics[width=\linewidth]{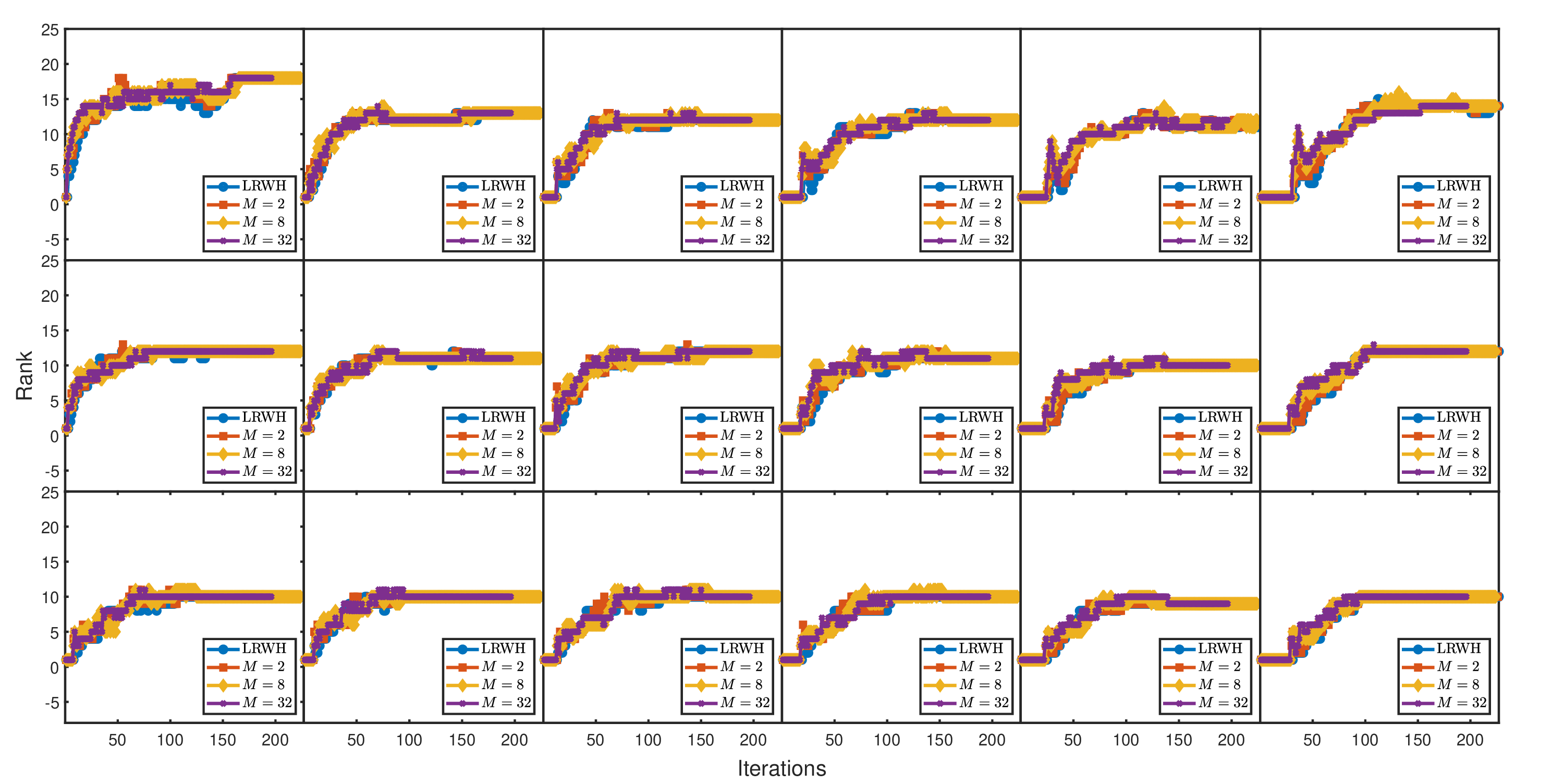}
\caption{Layered half-space problem. Distribution of ranks for the LRWH and LRAA methods in each of the $6\times 3$ blocks using the scheduling parameter $\theta=1$ and stopping tolerance $\epsilon^\star=10^{-3}$ for the memory parameters $M=2,8,16$ as well as the rank obtained from the full rank solution $W^n$ truncated at the tolerance at the $n$'th iteration. }
    \label{fig:ReflectingHalfspaceRanks}
\end{figure}

Layerings can also occur when the water-density varies with depth due to varying pressures and temperature. This is the physical mechanism behind the SOFAR channel \cite{Jenssen2011}. The SOFAR channel is a layer in the sea which is bounded along the top and bottom by layers with higher wave-speed, making it act as a waveguide. In this case, the waves generated in the channel will be trapped, leading to internal reflections. We model a simplified version of this setup as a layered half-space problem, mimicking reflections from the bounding layer. We set the wave speeds to be $c_1=0.7$ and $c_2=1$ as in the previous example.
\begin{figure}[]
    \centering
    \includegraphics[width=0.55\linewidth]{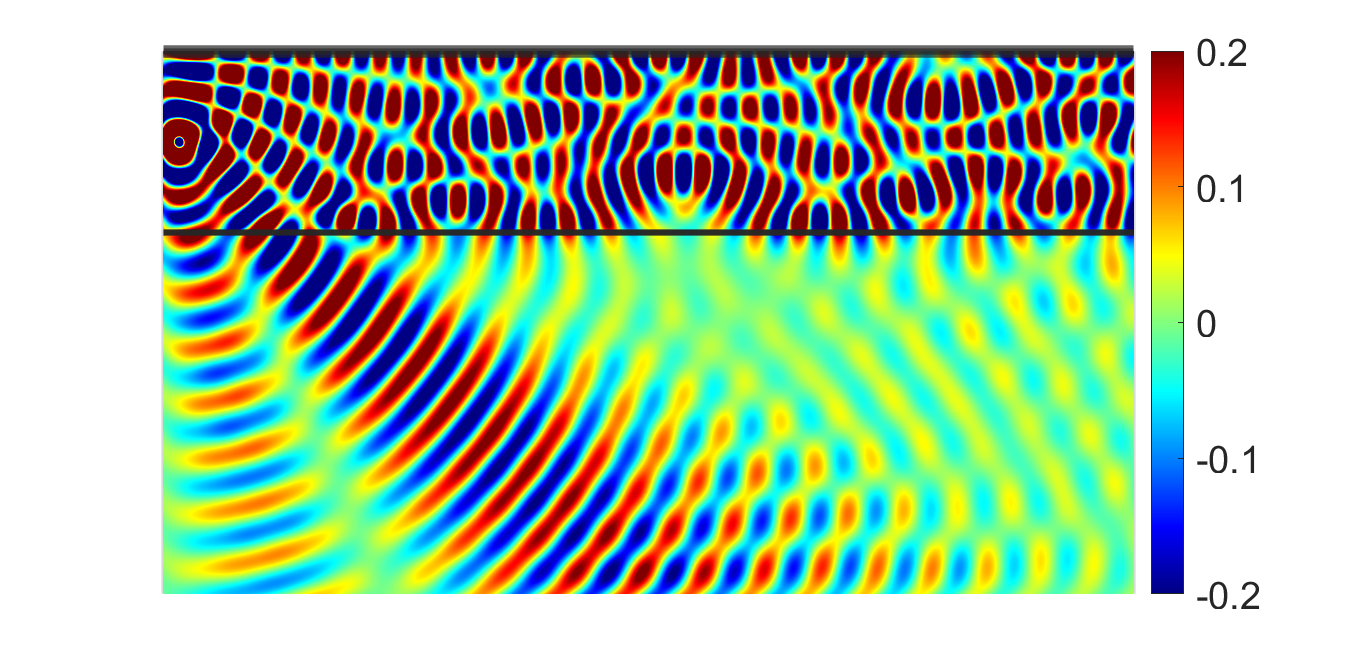}
   \includegraphics[width=0.4\linewidth]{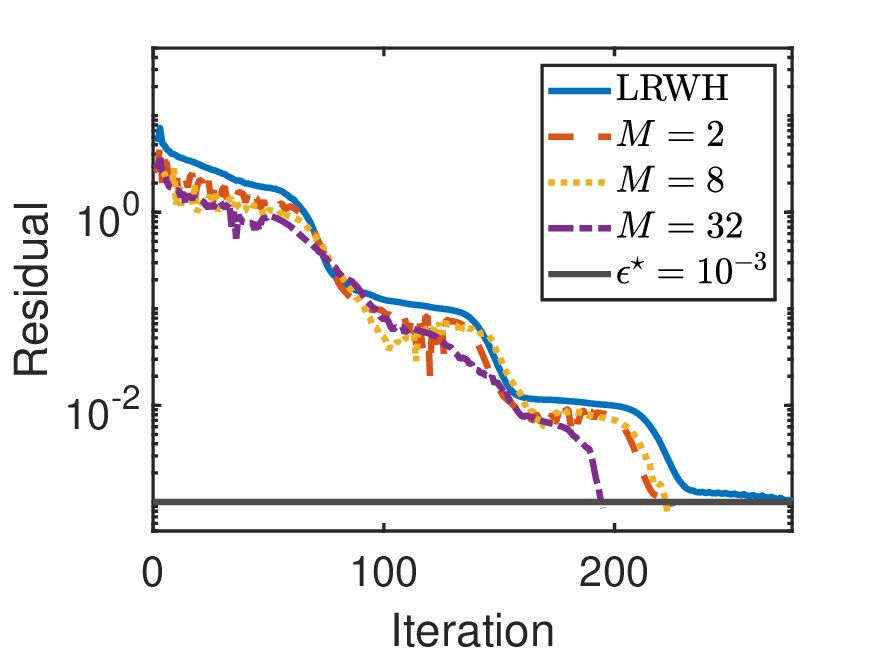}
    \caption{Layered half-space problem. The resulting Helmholtz solution in a discontinuous medium partitioned into $6\times 3$ blocks caused by a point source located above the interface (left) residual for the LRWH and $\operatorname{min}(\rho_G,\rho_X)$ obtained from LRAA using memory parameters $M=2,8,32$ with stopping tolerance $\epsilon^\star=10^{-3}$ and scheduling parameter $\theta=1$ (right).}
    \label{fig:ReflectingHalfspace}
\end{figure}
The resulting solution is shown together with the residuals in Figure~\ref{fig:ReflectingHalfspace} where the reflections in the top layer are clearly visible. We note that the LRWH has a relatively slow convergence compared with the layered free-space problem. We also note that the LRAA method indeed gives a clear acceleration. For $M=2$ and $M=8$ we have a gain of roughly 50 iterations, while increasing to \revone{$M=32$} results in a speedup of 80 iterations. The acceleration here is expected since we have trapped waves. The ranks are presented in Figure~\ref{fig:ReflectingHalfspaceRanks}. They follow the trend of having high rank close to the source as observed in the previous experiments. 
 
\subsection{Problems in 3D}
\label{subsec:3Dproblems}
We end the section by considering experiments in three dimensions, representing the solution using the tensor train format. The WaveHoltz iteration is generalized from two to three dimensions in a straightforward way using the rounding function in the \texttt{TT-toolbox}, we denote this method as the \textit{tensor train WaveHoltz} (TTWH) method. \revone{To our knowledge, there is no generalization of the fADI method to 3D tensors, so we choose a damping approach to model the open domain.}

\subsubsection{Timing and compression tests}
\begin{figure}[]
\centering
   \includegraphics[width=0.9\linewidth]{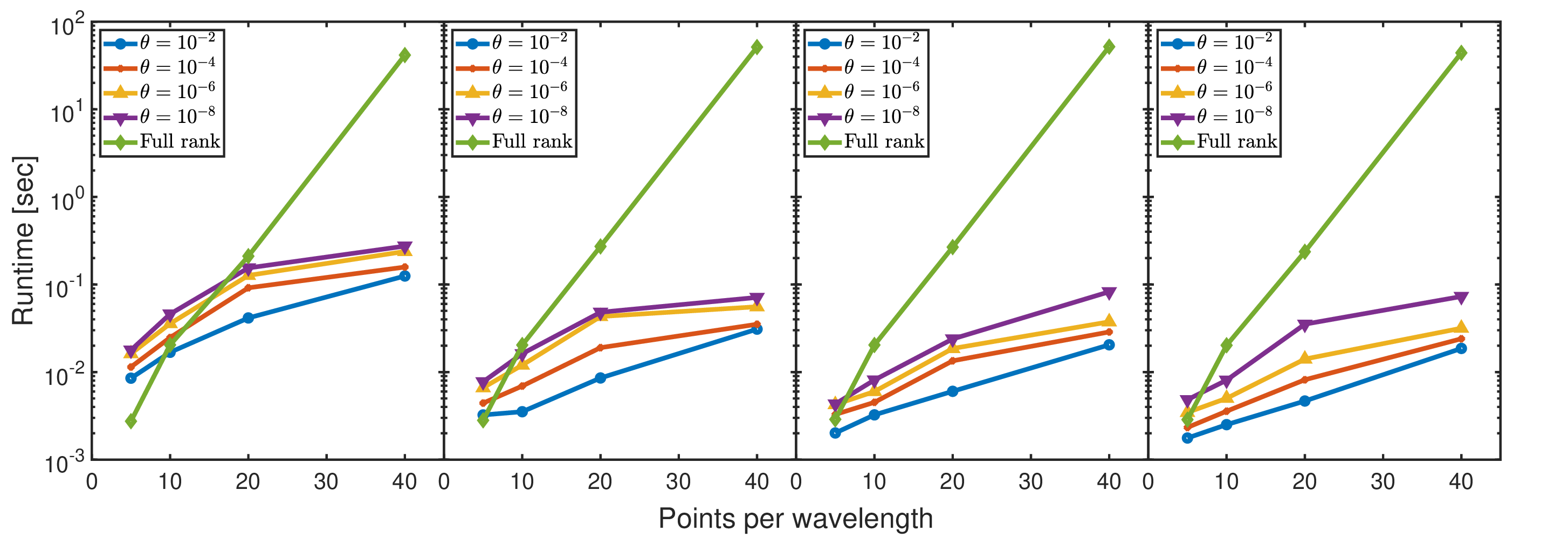}
    \includegraphics[width=0.9\linewidth]{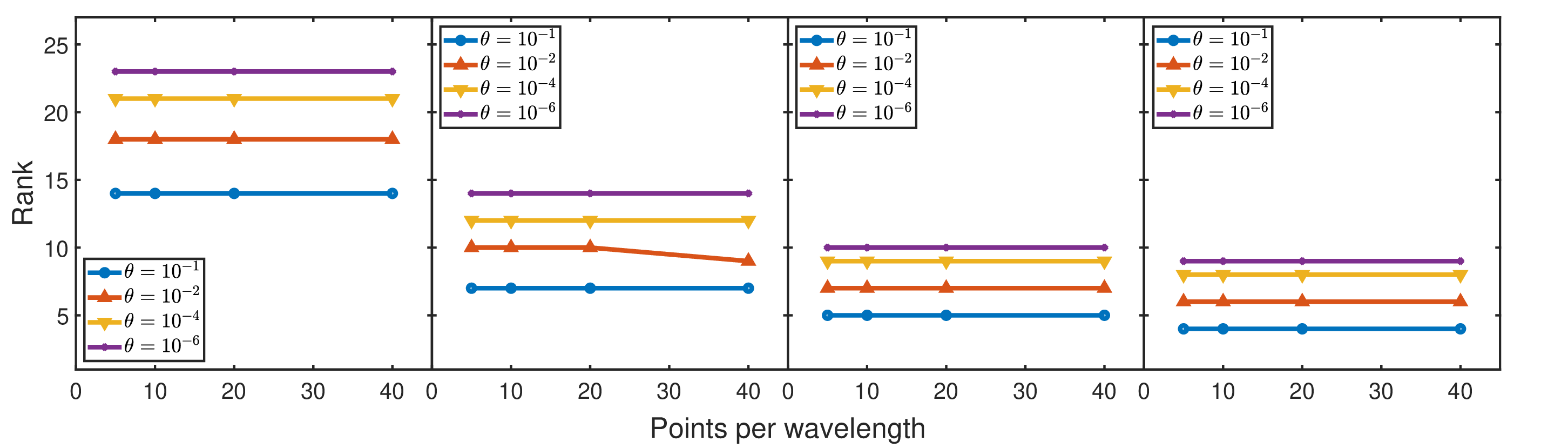}
\caption{Runtimes for taking a time step in each block (top)  maximal TT-ranks of the grid evaluation $\Gb$ (bottom) rounded by the blockwise tolerance $\epsilon=\sqrt{2}\theta/\|\Gb\|$ as functions of gridpoints per block.}
\label{fig:3DRuntimes}
\end{figure}

As in the two dimensional case, we begin by looking at the potential speedup by using a low-rank format compared to the full discretization. We begin by measuring the time it takes to apply the discrete Laplacian operator to a TT tensor $\Gb$ followed by the truncation operator for three different choices of tolerance levels. To make the rounding absolute and grid independent as in 2D we set $\delta=\sqrt{2}\epsilon/\|\Gb\|$ in \eqref{eq:TTrounding}. The runtime in the low-rank case is then compared with the time it would take to apply the discrete Laplacians to the vectorized solution. The domain is chosen as $\Omega=[0,4]\times[0,1]\times[0,1]$ and $\Gb$ the tensorized grid evaluation of the three dimensional Greens function of the open Helmholtz problem
\begin{equation}
    G(x,y,z)=\frac{1}{4\pi r(x,y,z)}e^{i\omega r(x,y,z)},\quad r(x,y,z)=\sqrt{(x-x_0)^2+(y-y_0)^2+(z-z_0)^2},
    \label{eq:3DGreen'sfnc}
\end{equation}
centered in the point $(x_0,y_0,z_0)=(-0.1,0.5,0.5)$. The runtimes for calculating $\Tc_\epsilon(\Gb+(\Delta t)^2[{\bf D}_x\Gb+{\bf D}_y\Gb+{\bf D}_z\Gb])$ in each block is shown in Figure \ref{fig:3DRuntimes} for different tolerances $\epsilon$ together with the runtime to perform the corresponding full rank calculation. We also present the maximal TT-ranks of $\Gb$ as functions of the number of points per wavelength used. Here we recall that ${\bf D}$ denotes the SBP operator $D_2$ converted to a TT-matrix. First, we note the large gain in runtime. Note that for $PPW=40$ the gain is substantial, almost two orders of magnitude compared to the full rank solver.
\subsubsection{Free space problem}
\begin{figure}[]
    \centering
    \includegraphics[width=0.45\linewidth]{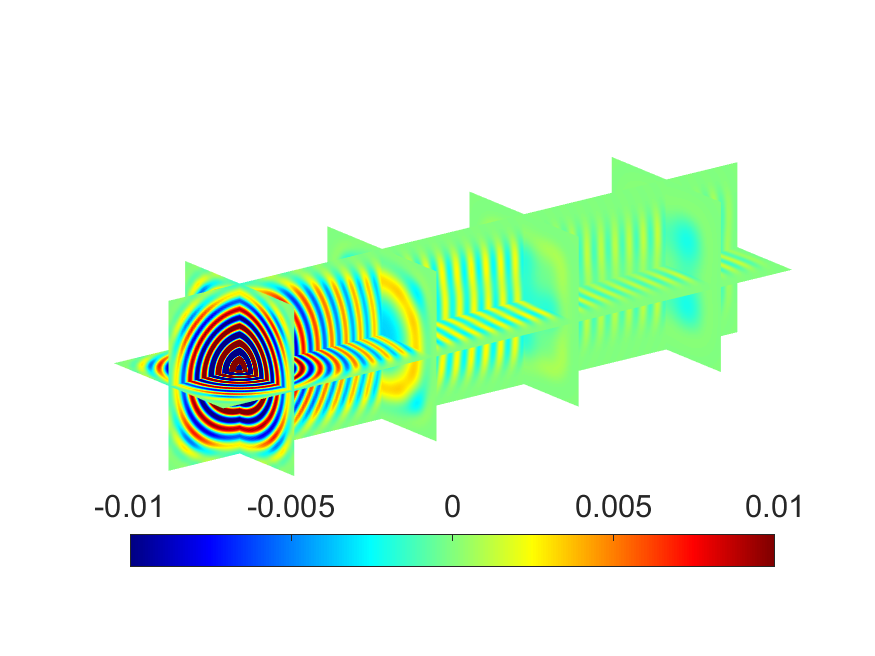}
    \includegraphics[width=0.45\linewidth]{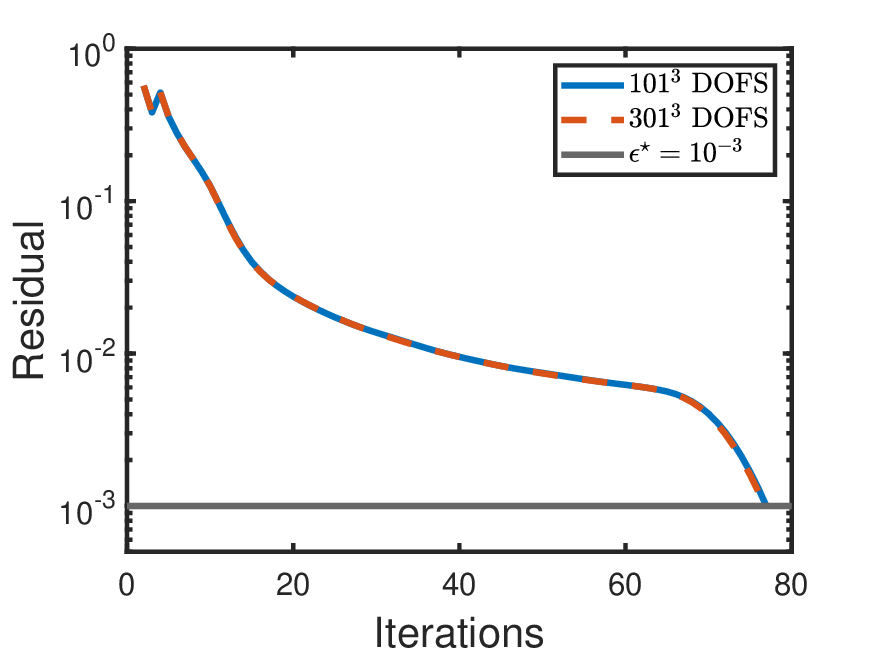}
    \caption{Free-space problem with unit wave speed solved using the TTWH method. The TTWH solution in $\Omega=[0,4]\times[0,1]\times[0,1]$ partitioned into $4\times 1\times 1$ blocks,  generated by a Gaussian point source located at $\xb_0=(0.5,0.5,0.5)$ with strong damping along all boundaries of the channel (left)  TTWH residuals obtained using the scheduling parameter $\theta=0.5$ with each block discretized by $101^3$ and $301^3$ gridpoints respectively and gridsize $h$. (right).}
    \label{fig:3DOutlowExSoln}
\end{figure}

In this experiment, we consider the TTWH method in the same domain $\Omega=[0,4]\times[0,1]\times[0,1]$ as in the previous experiment, but now with a Gaussian source located in the first block at the point $\xb_0=(1/4,1/2,1/2)$. To mimic an open domain we use a damping function rather than imposing outflow boundary conditions. The damping function is $\kappa(x,y,z)=\kappa_1(x)\kappa_2(y)\kappa_3(z)$ where
\begin{align*}
    \kappa_1(x)=50(e^{-100x^2}+e^{-100(x-4)^2}),\\
    \kappa_2(y)=50(e^{-100y^2}+e^{-100(y-1)^2}),\\
    \kappa_3(z)=50(e^{-100z^2}+e^{-100(z-1)^2}).
\end{align*}
We approximate the point source by the extension of a \eqref{eq:Gaussian} to three dimensions. 

We first consider blocks discretized using $101^3$ degrees of freedom, so that $\PPW=10$ and $h=10^{-2}$. We use $\rho^k_{ij}$ to determine the block-wise truncation tolerance. This is scaled by the norm of the current iterate and  grid size $h$ to obtain an absolute rounding. Moreover, we use scheduling to obtain the tolerance $\epsilon^{k+1}_{ij}=\max\{K,\sqrt{2}\theta h\rho^k_{ij}\}/\|\Wb_{ij}^k\|$. Note that we also use the maximum strategy here to prevent artificial rank-inflation when the residual is small. We set $K=10^{-4}$, $\theta=0.5$ and the stopping tolerance $\epsilon^\star=10^{-3}$, correpsonding to an integrated residual $h^{3/2}\epsilon^\star=10^{-6}$. The solution and residual shown in Figure~\ref{fig:3DOutlowExSoln}, demonstrating the expected low-rank behavior away from the point source in the form of simple wave fronts. We can also see that the TTWH iteration converges from the residual.

\begin{figure}[]
    \centering
   \includegraphics[width=0.8\linewidth]{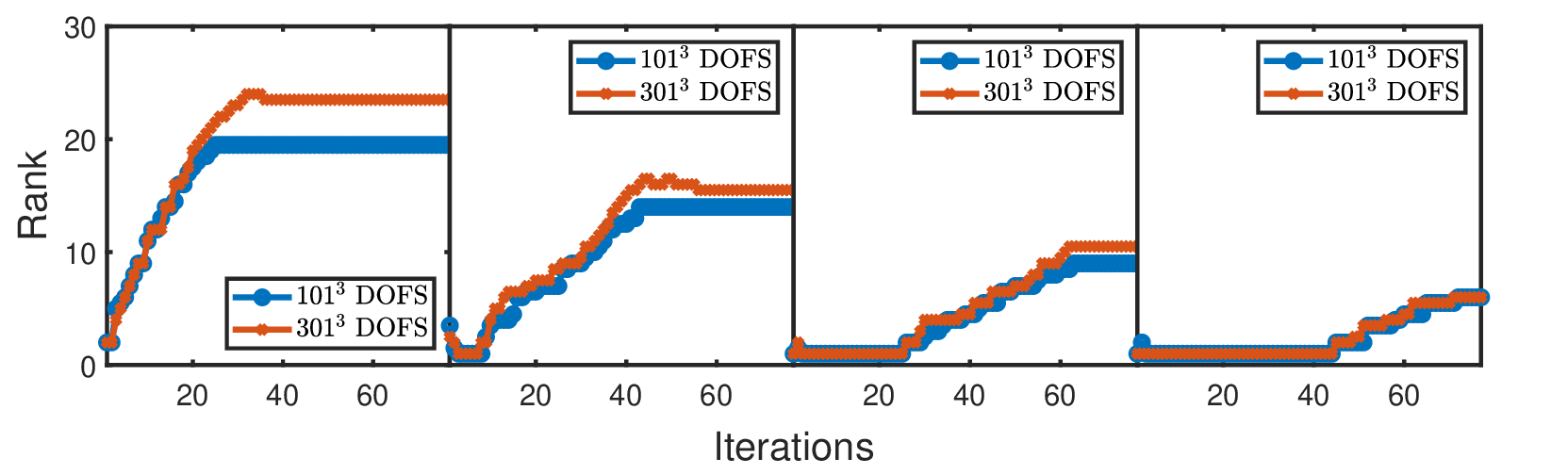}
    \caption{Free-space problem with unit wave speed solved using the TTWH method. Maximal TT ranks in each block discretizing each block using $101^3$ and $301^3$ degrees of freedom with scheduling parameter $\theta=0.5$.}
    \label{fig:3DOutflowExRanks}
\end{figure}

Since the rank of the solution is inherent to the underlying problem, if we fully resolve the solution the rank should remain the same even for finer meshes. We therefore consider the case when each block is discretized with $n=301^3=2.727\times 10^7$ gridpoints, which likely would be prohibitive using a standard full rank solver. We see from Figure~\ref{fig:3DOutflowExRanks} that the TT-ranks in each block remains similar for both discretizations, indicating a good resolution of the solution. The residual is also presented in Figure~\ref{fig:3DOutlowExSoln} and demonstrates a nearly identical convergence compared with the coarser discretization.
\subsubsection{Half-space problem}
\begin{figure}[]
    \centering
    \includegraphics[width=0.45\linewidth]{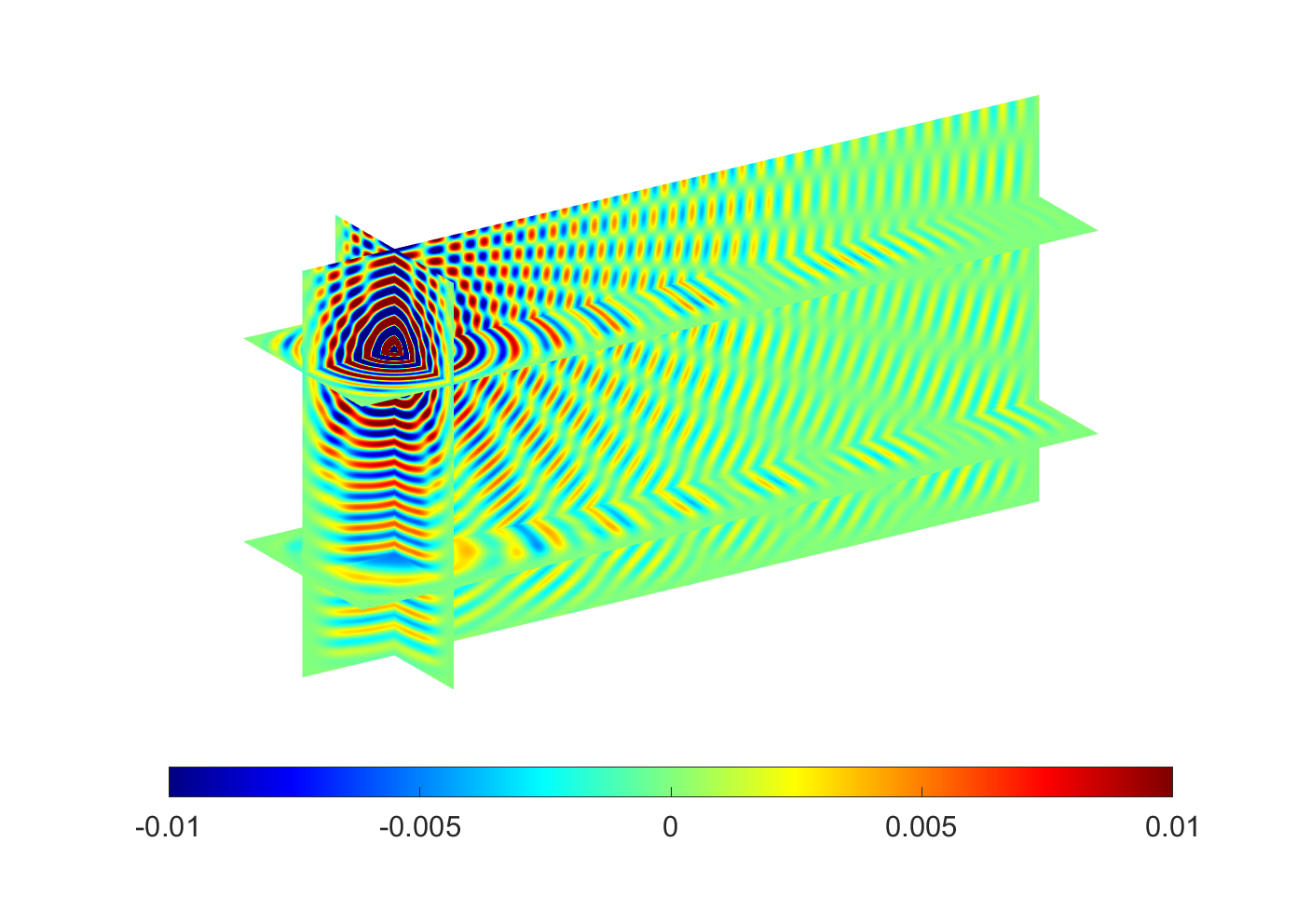}
    \includegraphics[width=0.45\linewidth]{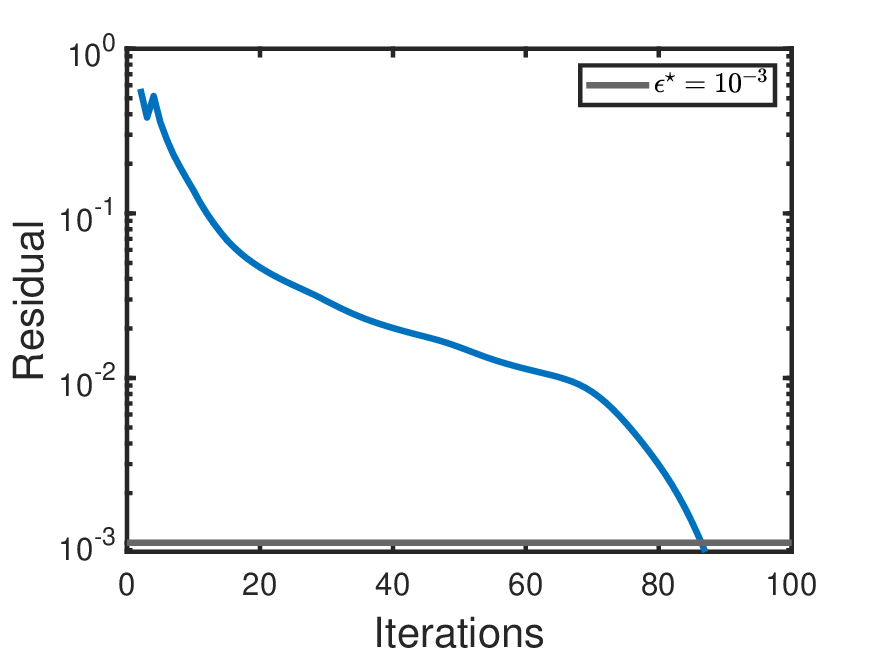}
    \caption{Mixed reflecting and outflow boundary conditions with unit wave speed solved using the TTWH method. The resulting Helmholtz solution as seen from the north-western corner of the computational domain. The domain is partitioned into $4\times 2\times 1$ blocks, each discretized with $101^3$ DOFS, generated by a point source below a reflecting water surface (left) TTWH residual using the scheduling parameter $\theta=0.5$ with each block discretized by $101^3$ gridpoints (right).}
    \label{fig:3DReflectingExSoln}
\end{figure}
\begin{figure}[ht!]
    \centering
 \includegraphics[width=0.8\linewidth]{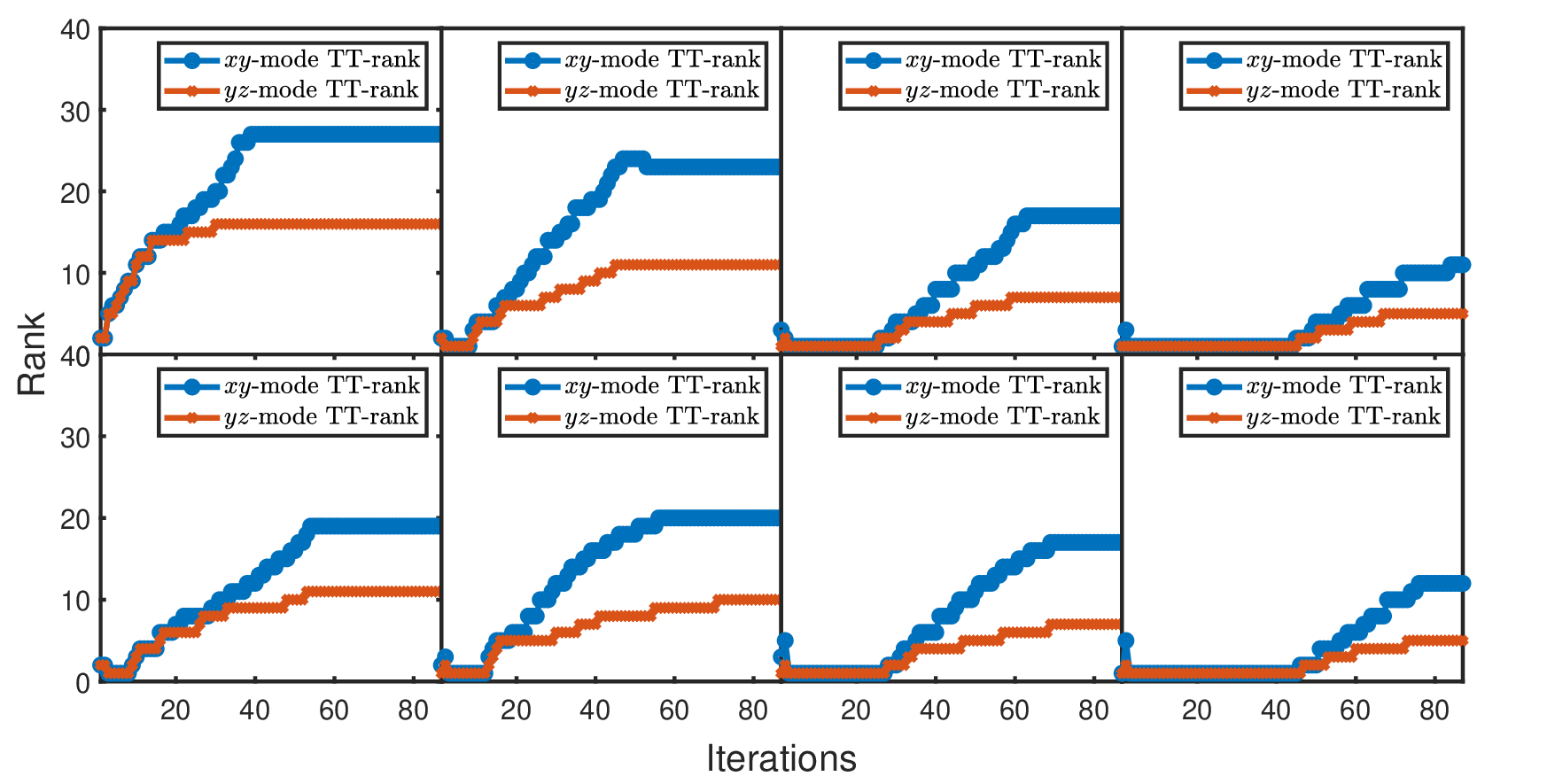}
    \caption{Half-space problem with unit wave speed solved using the TTWH method. TT ranks of the $xy$ and $yz$ modes respectively in each of the $4\times 2\times 1$ blocks, each discretized with $101^3$ DOFS, using the scheduling parameter $\theta=0.5$.}
    \label{fig:3DReflectingRanks}
\end{figure}

As in the two dimensional setting, we can model interations with the water surface by imposing a reflecting boundary condition on the top surface of the domain. We do so by simply removing the damping function along the northern boundary. We partition the computational domain into $4\times 2\times 1$ blocks, each discretized by $101^3$ gridpoints \revone{so that $h=10^{-2}$, maintaining $\text{PPW}=30$. Finally, we use the same truncation strategy as in the free-space problem.} The resulting solution and residual is presented in Figure~\ref{fig:3DReflectingExSoln}. We see that the solution behaves as expected far away from the point source, while the solution along the leftmost boundary is highly damped due to the short distance to the pointsource. We also note that the convergence is much slower than in the case of a fully open domain, which is consistent with the results in two dimensions. Finally, we present the rank distribution throughout each block in Figure~\ref{fig:3DReflectingRanks}. Also in this case we see that the rank grows almost monotonically, with the exception of an increase during the initial iterations when the solution barely has entered the blocks. The difference in ranks between xy- and yz- modes is likely due to the difference in structural complexity of the solution in the xy and yz planes, as seen in Figure~\ref{fig:3DReflectingExSoln}.

\subsubsection{Effect of distance to source on compression}
\begin{figure}
    \centering
    \includegraphics[width=0.8\linewidth]{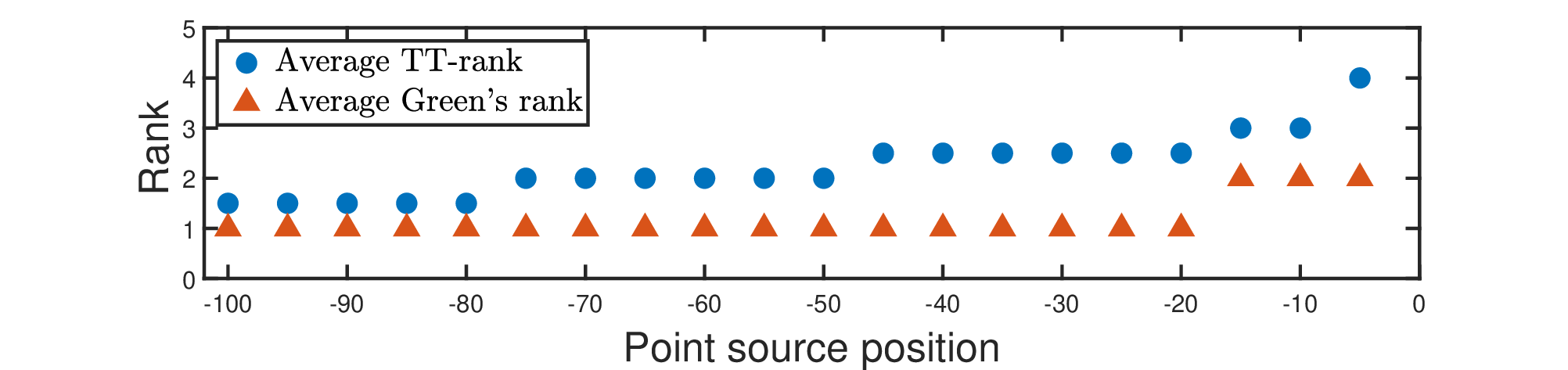}
    \caption{Moving point source. Average TT-ranks of the numerical solution and Green's function in $[0,1]^3$ with damping, discretized by $101^3$ gridpoints as a function of the point source position along the $x$-axis using the scheduling parameter $\theta=0.5$}.
    \label{fig:MovingCubeRanks}
\end{figure}
In the final experiment we consider the effect of truncation in a long domain. To model this, we consider a unit cube $[0,1]^3$ as if it was an element in a longer domain by imposing a Dirichlet boundary condition on the western face of the cube and place the source outside the cube. The damping function is chosen to localize to the boundaries of the cube. We set the Dirichlet data as the Green's function \eqref{eq:3DGreen'sfnc}, moving the source further and further away by varying its position along the x-axis then models the cube being further and further downstream in the channel. We follow the scheduling strategy used in the previous experiments. Since the computational cost to evaluate the TT SVD is $\mathcal{O}(3nr^3)$ compared with the full rank cost of $\mathcal{O}(n^3)$, we see that the low rank method is more efficient $r\sim (n^2/3)^{1/3}$. Setting \revone{$n=101$ so that $h=10^{-2}$ and $\PPW=10$}, we get that there is a tradeoff when $r=15$. We present the final TT-ranks in Figure~\ref{fig:MovingCubeRanks} as a function of x-coordinate of the point source, together with the rank of the analytical Green's function \eqref{eq:3DGreen'sfnc} truncated with the final tolerance used in the TTWH solver. It is clear that the maximal TT-rank is always bounded by the tradeoff rank, indicating that the method always remains more efficient than a full rank solver for this problem setup. Moreover, we see that the TTWH ranks approach the analytical ranks.
\section{Conclusion}\label{sec:conclusion}
In this work, we have developed LRWH, a low-rank WaveHoltz method, for simulating underwater acoustics. The solver employs a low-rank representation of the solution, using singular value decomposition or tensor trains in two and three dimensions respectively. To control rank growth we apply an explicit step-truncation strategy in combination with scheduling of the truncation level using the residual. 

We find that in 3D the algorithm gives large compression compared to a full rank solver, allowing for simulation of larger problems. In 2D, the benefits of using the low-rank format is limited. Moreover, for the open problems considered here, the WaveHoltz iteration itself is already highly efficient and therefore acceleration by Anderson does not improve iteration count significantly. 

Several numerical examples were presented, showing that the method can handle various configurations of material properties and boundary conditions. An important property of the method is that throughout the iteration the rank is monotonically increasing towards the final rank. Very little, if any, intermediate rank inflation is observed.

\revonesecondgo{To improve efficiency, in future work we will consider hybridizing full- and low-rank methods to exploit the decrease in rank as a function of distance to source. Another important extension is to curvilinear grids, being able to consider realistic bathymetry in underwater acoustics applications. The main challenge in this setting is that the spatially varying coefficients in the coordinate transformation couple neighboring grid values in a way that destroys the separable structure used in the SVD-based rank truncation. The solution is to replace the truncation step with a cross-approximation technique \cite{Appelo2025}, which can approximate the action of the operator along subsets of rows and columns to recover a low-rank representation without requiring separability.} 
\section*{Acknowledgements}
Part of the computations were enabled by resources provided by the National Academic Infrastructure for Supercomputing in Sweden (NAISS) and HPC2N (project ID: hpc2n2025-198), partially funded by the Swedish Research Council through grant agreement no. 2022-06725.

Granath acknowledges support from the foundations managed by The Royal Swedish Academy of Sciences, project number MA2024-0086, and the Kempe foundation.

Appel{\"o} is supported in part by U.S. Department of Energy, Oﬃce of Science,
Advanced Scientific Computing Research (ASCR), under Award Number DE-SC0025424, by the National Science Foundation under grant NSF DMS-2208164, NSF DMS-2436319, and Virginia Tech. This material is based upon work supported by the National Science Foundation under Grant No. DMS-1928930 while Appel\"{o} was in residence at the Simons Laufer Mathematical Sciences
Institute in Berkeley, California, during the Fall 2025 semester.

\bibliographystyle{elsarticle-num} 
\bibliography{refs}

@article{ROTEM2026,
title = {Convergence of the semi-discrete {WaveHoltz} iteration},
journal = {Journal of Computational Physics},
volume = 558,
pages = 114882,
year = 2026,
issn = {0021-9991},
doi = {https://doi.org/10.1016/j.jcp.2026.114882},
author = {Amit Rotem and Olof Runborg and Daniel Appel\"{o}}
}

@PhdThesis{sabino2006,
  title={Solution of large-scale {L}yapunov equations via the block modified {S}mith method},
  author={Sabino, John},
  school={Rice University, Houston},
  year=2006
}

@article{LowRankFE,
	author = {Bebendorf, Mario and Hackbusch, Wolfgang},
	journal = {Numerische Mathematik},
	number = {1},
	pages = {1--28},
	title = {Existence of $\mathcal{H}$-matrix approximants to the inverse {FE}-matrix of elliptic operators with $L_\infty$-coefficients},
	volume = {95},
	year = {2003}}

@article{parabolic_sep,
author = {Engquist, Bj{\"o}rn and Ying, Lexing},
title = {Fast Directional Multilevel Algorithms for Oscillatory Kernels},
journal = {SIAM Journal on Scientific Computing},
volume = 29,
number = 4,
pages = {1710-1737},
year = 2007}

@article{no_low_rank_Engq,
	Author = {B. Engquist and H. Zhao},
	Doi = {10.1002/cpa.21755},
	Journal = {Communications on Pure and Applied Mathematics},
	Number = 11,
	Pages = {2220-2274},
	Title = {Approximate Separability of the {G}reen's Function of the {H}elmholtz Equation in the High Frequency Limit},
	Volume = 71,
	Year = 2018}

@article{Almquist2019,
  title={Order-preserving interpolation for summation-by-parts operators at nonconforming grid interfaces},
  author={Almquist, Martin and Wang, Siyang and Werpers, Jonatan},
  journal={SIAM Journal on Scientific Computing},
  volume={41},
  number={2},
  pages={A1201--A1227},
  year={2019},
  publisher={SIAM}
}

@article{Appelo2019,
  title={An energy-based discontinuous {G}alerkin method for coupled elasto-acoustic wave equations in second-order form},
  author={Appel\"{o}, Daniel and Wang, Siyang},
  journal={International Journal for Numerical Methods in Engineering},
  volume={119},
  number={7},
  pages={618--638},
  year={2019},
  publisher={Wiley Online Library}
}

@article{Appelo2020elastic,
  title={{E}l-{W}ave{H}oltz: A time-domain iterative solver for time-harmonic elastic waves},
  author={Appel\"{o}, Daniel and Garcia, Fortino and Loya, Allen Alvarez and Runborg, Olof},
  journal={Computer Methods in Applied Mechanics and Engineering},
  volume={401},
  pages={115603},
  year={2022},
  publisher={Elsevier}
}

@article{Appelo2020,
  title={WaveHoltz: Iterative solution of the {H}elmholtz equation via the wave equation},
  author={Appel\"{o}, Daniel and Garcia, Fortino and Runborg, Olof},
  journal={SIAM Journal on Scientific Computing},
  volume={42},
  number={4},
  pages={A1950--A1983},
  year={2020},
  publisher={SIAM}
}

@article{Appelo2025,
      title={{lrAA}: {L}ow-{R}ank {A}nderson {A}cceleration}, 
      author={Daniel Appel\"{o} and Yingda Cheng},
      year=2026,
      journal={SISC (accepted)} 
}

@article{Appelo2025Rule,
  title={A Rule of Thumb for Choosing Points-Per-Wavelength for Finite Difference Approximations of Helmholtz Problems},
  author={Appel{\"o}, Daniel and Banks, Jeffrey W and Henshaw, William D and Schwendeman, Donald W},
  journal={Journal of Computational Physics},
  pages={114703},
  year={2026},
  publisher={Elsevier}
}

@article{Bachmayr2023,
  title={Low-rank tensor methods for partial differential equations},
  author={Bachmayr, Markus},
  journal={Acta Numerica},
  volume={32},
  pages={1--121},
  year={2023},
  publisher={Cambridge University Press}
}

@article{Benner2009,
  title={On the ADI method for Sylvester equations},
  author={Benner, Peter and Li, Ren-Cang and Truhar, Ninoslav},
  journal={Journal of Computational and Applied Mathematics},
  volume={233},
  number={4},
  pages={1035--1045},
  year={2009},
  publisher={Elsevier}
}

@article{Carpenter1994,
  title={Time-stable boundary conditions for finite-difference schemes solving hyperbolic systems: methodology and application to high-order compact schemes},
  author={Carpenter, Mark H and Gottlieb, David and Abarbanel, Saul},
  journal={Journal of Computational Physics},
  volume={111},
  number={2},
  pages={220--236},
  year={1994},
  publisher={Elsevier}
}

@article{Dektor2021,
  title={Rank-adaptive tensor methods for high-dimensional nonlinear PDEs},
  author={Dektor, Alec and Rodgers, Abram and Venturi, Daniele},
  journal={Journal of Scientific Computing},
  volume={88},
  number={2},
  pages={36},
  year={2021},
  publisher={Springer}
}

@article{Ernst2011,
  title={Why it is difficult to solve {H}elmholtz problems with classical iterative methods},
  author={Ernst, Oliver G and Gander, Martin J},
  journal={Numerical analysis of multiscale problems},
  pages={325--363},
  year={2011},
  publisher={Springer}
}

@article{Erlangga2008,
  title={Advances in iterative methods and preconditioners for the {H}elmholtz equation},
  author={Erlangga, Yogi A},
  journal={Archives of Computational Methods in Engineering},
  volume={15},
  pages={37--66},
  year={2008},
  publisher={Springer}
}

@article{Hagstrom1999,
  title={Radiation boundary conditions for the numerical simulation of waves},
  author={Hagstrom, Thomas},
  journal={Acta numerica},
  volume={8},
  pages={47--106},
  year={1999},
  publisher={Cambridge University Press}
}

@article{Hagstrom2012,
  title={Grid stabilization of high-order one-sided differencing II: Second-order wave equations},
  author={Hagstrom, Thomas and Hagstrom, George},
  journal={Journal of Computational Physics},
  volume={231},
  number={23},
  pages={7907--7931},
  year={2012},
  publisher={Elsevier}
}

@article{Hicken2013,
  title={Summation-by-parts operators and high-order quadrature},
  author={Hicken, Jason E and Zingg, David W},
  journal={Journal of Computational and Applied Mathematics},
  volume={237},
  number={1},
  pages={111--125},
  year={2013},
  publisher={Elsevier}
}

@book{Hovem2011,
  title={Ray trace modeling of underwater sound propagation. Documentation and use of the PlaneRay model},
  author={Hovem, Jens Martin},
publisher={IntechOpen},
  year={2011}
}

@book{Jenssen2011,
  title={Computational ocean acoustics},
  author={Jensen, Finn B and Kuperman, William A and Porter, Michael B and Schmidt, Henrik and Tolstoy, Alexandra},
  volume={2011},
  year={2011},
  publisher={Springer}
}

@article{Koch2007dynamical,
  title={Dynamical low-rank approximation},
  author={Koch, Othmar and Lubich, Christian},
  journal={SIAM Journal on Matrix Analysis and Applications},
  volume={29},
  number={2},
  pages={434--454},
  year={2007},
  publisher={SIAM}
}

@article{Kreiss1972,
  title={Comparison of accurate methods for the integration of hyperbolic equations},
  author={Kreiss, Heinz-Otto and Oliger, Joseph},
  journal={Tellus},
  volume={24},
  number={3},
  pages={199--215},
  year={1972},
  publisher={Taylor \& Francis}
}

@incollection{Kreiss1974,
  title={Finite element and finite difference methods for hyperbolic partial differential equations},
  author={Kreiss, H-O and Scherer, Godela},
  booktitle={Mathematical aspects of finite elements in partial differential equations},
  pages={195--212},
  year={1974},
  publisher={Elsevier}
}

@article{Mattsson2004,
  title={Summation by parts operators for finite difference approximations of second derivatives},
  author={Mattsson, Ken and Nordstr{\"o}m, Jan},
  journal={Journal of Computational Physics},
  volume={199},
  number={2},
  pages={503--540},
  year={2004},
  publisher={Elsevier}
}

@article{Mattsson2009,
  title={Stable boundary treatment for the wave equation on second-order form},
  author={Mattsson, Ken and Ham, Frank and Iaccarino, Gianluca},
  journal={Journal of Scientific Computing},
  volume={41},
  pages={366--383},
  year={2009},
  publisher={Springer}
}

@article{Mattsson2008,
  title={Stable and accurate wave-propagation in discontinuous media},
  author={Mattsson, Ken and Ham, Frank and Iaccarino, Gianluca},
  journal={Journal of Computational Physics},
  volume={227},
  number={19},
  pages={8753--8767},
  year={2008},
  publisher={Elsevier}
}

@article{Oseledets2011,
  title={Tensor-train decomposition},
  author={Oseledets, Ivan V},
  journal={SIAM Journal on Scientific Computing},
  volume={33},
  number={5},
  pages={2295--2317},
  year={2011},
  publisher={SIAM}
}

@article{Pekeris1948,
  title={Theory of propagation of explosive sound in shallow water},
  author={Pekeris, Chaim Leib},
  year={1948}
}

@article{Peng2022,
  title={{EM}-{W}ave{H}oltz: A flexible frequency-domain method built from time-domain solvers},
  author={Peng, Zhichao and Appel{\"o}, Daniel},
  journal={IEEE Transactions on Antennas and Propagation},
  volume={70},
  number={7},
  pages={5659--5671},
  year={2022},
  publisher={IEEE}
}

@article{Petersen2008,
  title={The matrix cookbook},
  author={Petersen, Kaare Brandt and Pedersen, Michael Syskind and others},
  journal={Technical University of Denmark},
  volume={7},
  number={15},
  pages={510},
  year={2008}
}

@article{Tappert2005,
  title={The parabolic approximation method},
  author={Tappert, Fred D},
  journal={Wave propagation and underwater acoustics},
  pages={224--287},
  year={2005},
  publisher={Springer}
}

@article{Yang2022,
  title={{A}nderson acceleration based on the $\mathcal{H}^{-s}$ Sobolev norm for contractive and noncontractive fixed-point operators},
  author={Yang, Yunan and Townsend, Alex and Appel\"{o}, Daniel},
  journal={Journal of Computational and Applied Mathematics},
  volume={403},
  pages={113844},
  year={2022},
  publisher={Elsevier}
}
\appendix
\section{Enforcing time dependent boundary conditions}
\label{subsec:fADI}
In this section, we provide a short introduction to the method used to impose the outflow conditions in matrix form. We approximate the time derivatives in the boundary conditions of \eqref{eq:fullrankdiscU}-\eqref{eq:fullrankdiscV} with second order centered differences using a time step $\Delta t$ and let $t^k=k\Delta t$. Then, the fully discrete equation for $\Umf$ can be written as
\begin{equation}
    \Umf^{k+1}+\At\Umf^{k+1}+\Umf^{k+1}\Bt^T=R,
    \label{eq:FRSylvester}
\end{equation}
where $R$ depends on $\Umf$ and $\Vmf$ at times $t^k$ and $t^{k-1}$ and the matrices $\At$ and $\Bt$ depends on the current block we are considering. A key difference of this scheme compared with the vector-based solvers is that we cannot obtain $\Umf^{k+1}$ by a direct solve. Instead, $\Umf^{k+1}$ is given in terms of a matrix equation. There is an extensive literature on solving these types of equations, and we resort to using the \textit{factored alternating direction iteration} (fADI) method \cite{Benner2009}. This method is particularly appealing as it is tailored towards equations where the matrices $A$ and $B$ have low-rank. To reformulate the equation for the fADI method, we rewrite the left hand side of \eqref{eq:FRSylvester} as
\begin{equation}
   \Umf^{k+1}+\At\Umf^{k+1}+\Umf^{k+1}\Bt^T=(\frac{1}{2}I_n+\At)\Umf^{k+1}-\Umf^{k+1}(-\frac{1}{2}I_n-\Bt)^T=A \Umf^{k+1}-B \Umf^{k+1}=GF^T,
   \label{eq:FRSylvester2}
\end{equation}
where $I_n\in\Rbb^{n\times n}$ denotes the identity matrix. Since this method will be used in a low-rank framework, we assume that the matrix $R$ can be factored as $R=GF^T$ for some matrices $F,G\in\Rbb^{n\times r}$ where $r\leq n$. Then, given two sets of \textit{shift parameters} $\{\alpha_i\}_{i=1}^n,\{\beta_i\}_{i=1}^n$ the fADI iteration can be written as below in Algorithm~\ref{algorithm:fADI}. Note that we use the MATLAB notation $Z(:,i)$ to denote the $i$'th column and $Z(j,:)$ the $j$'th row. Moreover, the shift parameters can in principle be chosen arbitrarily, but it can be shown that they greatly affect the speed of convergence of the method \cite{Benner2009}. In particular, it is shown in \cite{sabino2006} that the coefficients can be chosen to solve the minimization problem
\begin{equation}
    \min_{\substack{\alpha_i\in\Cbb\\\beta_i\in\Cbb}}\max_{\substack{\lambda_i\in\operatorname{eig}(\At)\\\lambda'_i\in\operatorname{eig}(\Bt)}}\prod_{i=1}^n\bigg|\frac{(\alpha_i-\lambda_i)(\alpha_i-\lambda'_i)}{(\alpha_i-\lambda'_i)(\beta_i-\lambda_i)}\bigg|,
    \label{eq:ADIerror}
\end{equation}
\begin{algorithm}[!ht]
 \DontPrintSemicolon
    \caption{fADI for Sylvester Equation $A\Umf^{k+1}-\Umf^{k+1}B^T=GF^T$}

    \label{algorithm:fADI}
 
    \KwIn{$A,B\in\Rbb^{n\times n}$, $F,G\in\Rbb^{n\times r}$, ADI shifts $\{\alpha_i\},\{\beta_i\}$, number of fADI iterations $m$}
    \KwOut{Factors $Z,Y\in\Rbb^{m\times mr}$ and $D\in\Rbb^{mr\times mr}$ such that $\Umf^{k+1}=ZDY^T$ approximately solves $A\Umf^{k+1}-\Umf^{k+1}B^T=GF^T$}
    $Z(:,1:r)=(A-\beta_1I_n)^{-1}G$, $(Y^T)(:,1:r)=F^T(B-\alpha_1I_n)^{-1}$\\
    \For{$i=1,2,...,m$}{
    $Z(:,ir+1:(i+1)r)=Z(:,(i-1)r+1:ir)+(\beta_{i+1}-\alpha_i)(A-\beta_{i+1}I_n)^{-1}Z(:,(i-1)r+1:ir)$\\
    $(Y^T)(ir+1:(i+1)r,:)=(Y^T)((i-1)r+1:ir,:)+(\alpha_{i+1}-\beta_i)(Y^T)((i-1)r+1:ir,:))(B-\alpha_{i+1}I_n)^{-1}$\\    
    }
    $D=\operatorname{diag}((\beta_1-\alpha_1)I_r,\cdots,(\beta_k-\alpha_k)I_r)$\\
    Return $\Umf^{k+1}=ZDY^T$.
\end{algorithm}
where $\operatorname{eig}(\cdot)$ denotes the spectrum of a matrix. For our problem, we know $A=\ab\ab^T+\bb\bb^T$ and $B=\cb\cb^T+\db\db^T$ for some vectors $\ab,\bb,\cb,\db\in\Rbb^n$ satisfying $\ab^T\bb=\cb^T\db=0$. Since $\ab\ab^T$ and $\bb\bb^T$ are rank one matrices with $\ab^T\bb=0$, they will each contribute one unique eigenvalue to $\At$, which therefore will have the unique eigenvalues $\lambda_1=1/2+\ab^T\ab, \lambda_2=1/2+\bb^T\bb$ and then $\lambda_i=1/2$ with multiplicity $n-2$. Similarly we get the eigenvalues $\lambda'_1=-1/2-\cb^T\cb,\lambda'_2=-1/2-\db^T\db$ and $\lambda'_i=-1/2$ for $\Bt$. Choosing the shift parameters equal to the eigenvalues will make the iterations converge in three iterations. We can therefore represent the solution by analytic formula provided in the following lemma.

\begin{lemma}
    \label{lemma:fADIlemma} Let $\ab,\bb,\cb,\db\in\Rbb^n$ satisfy $\ab^T\bb=\cb^T\db=0$, $R\in\Rbb^{n\times n}$ a matrix with the SVD representation $USV^T$, $U,V\in\Rbb^{n\times r}$, $S\in\Rbb^{r\times r}$ and $\epsilon>0$ the truncation parameter. Then, the solution $X$ to the Sylvester equation
    $$AX-XB^T=R$$
    where $A=\frac{1}{2}I_n+\ab\ab^T+\bb\bb^T$ and $B=-\frac{1}{2}-\cb\cb^T-\db\db^T$ can be represented in terms of a low-rank sum
    $$X=\Tc_\epsilon^{sum}\bigg(\sum_{j=1}^3Z_i(\beta_i-\alpha_i)I_{r\times r}Y_i^T\bigg),$$
    where $\{\alpha_i\}_{i=1}^3,\{\beta_i\}_{i=1}^3$ denote the set of unique eigenvalues of $A$ and $B$. Moreover, the factor matrices $Z_i$ and $Y_i$ are given by

    \begin{align*}
        &Z_1=\begin{bmatrix}
            U &\ab &\bb
        \end{bmatrix}\begin{bmatrix}   
        \frac{\sqrt{S}}{1+\cb^T\cb} & \gamb_a^T\sqrt{S} &\gamb_b^T\sqrt{S} \end{bmatrix}\\
       &Y_1=\begin{bmatrix}
           V & \cb & \db
       \end{bmatrix}\begin{bmatrix}
           \frac{\sqrt{S}}{1+\ab^T\ab} & \gamb_c^T\sqrt{S} &\gamb_d^T\sqrt{S}
       \end{bmatrix}      
    \end{align*}
    and
    $Z_2=(\kz I_n+\kzz \ab\ab^T+\kzzz\bb\bb^T)Z_1$, $Z_3=(\kzt I_n+\kzzt\ab\ab^T+\kzzzt\bb\bb^T)Z_1$, $Y_2=(\ky I_n+\kyy\cb\cb^T+\kyyy\db\db^T)Y_1$ and $Y_3=(\kyt I_n+\kyyt\cb\cb^T+\kyyyt\db\db^T)Y_1$, where the coefficients are given by
    \begin{alignat*}{3}
    & \kz=-\frac{\ab^T\ab}{1+\db^T\db} \quad&& \kzz=\frac{1}{1+\db^T\db} \quad&& \kzzz=\frac{1+\ab^T\ab+\db^T\db}{(1+\db^T\db)(1+\bb^T\bb+\db^T\db)},\\ 
    & \kzt=\frac{(\ab^T\ab)(\bb^T\bb)}{1+\db^T\db}\quad && \kzzt=-\frac{\bb^T\bb}{1+\db^T\db}\quad && \kzzzt=-\frac{\ab^T\ab}{1+\db^T\db},\\
    & \ky=-\frac{\cb^T\cb}{1+\bb^T\bb}\quad && \kyy=\frac{1}{1+\bb^T\bb}\quad && \kyyy=\frac{1+\bb^T\bb+\cb^T\cb}{(1+\bb^T\bb)(1+\bb^T\bb+\db^T\db)},\\
    & \kyt=\frac{(\cb^T\cb)(\db^T\db)}{(1+\bb^T\bb)}\quad &&\kyyt=-\frac{\db^T\db}{1+\bb^T\bb}\quad && \kyyyt=-\frac{\cb^T\cb}{1+\bb^T\bb},
    \end{alignat*}
    and the vectors $\gamb_a,\gamb_b,\gamb_c$ and $\gamb_d$ by 
    \begin{alignat*}{2}
    & \gamb_a=-\frac{\ab^TU}{(1+\cb^T\cb)(1+\ab^T\ab+\cb^T\cb)} \quad&& \gamb_b=-\frac{\bb^TU}{(1+\cb^T\cb)(1+\bb^T\bb+\cb^T\cb)},\\
    &\gamb_c=-\frac{\cb^TV}{(1+\ab^T\ab)(1+\ab^T\ab+\cb^T\cb)} \quad&& \gamb_d=-\frac{\db^TV}{(1+\ab^T\ab)(1+\ab^T\ab+\db^T\db)}.       
    \end{alignat*}

\end{lemma}

\section{Proof of Lemma~\ref{lemma:fADIlemma}}
\label{proof:LRFADI}
We derive the exact solutions to the equation
$$AX-XB^T=USV^T,$$
where $USV^T$ is given $U,V\in\Rbb^{n\times r}$ and $S\in\Rbb^r$, $A=\frac{1}{2}I_n+\ab\ab^T+\bb\bb^T$ and $B=-\frac{1}{2}I_n-\cb\cb^T-\db\db^T$. The vectors $\ab,\bb,\cb$ and $\db$ satisfy $\ab^T\bb=\cb^T\db=0$. Recall that the eigenvalues of $A$ are given by $\lambda_1=\frac{1}{2}+\ab^T\ab,\lambda_2=\frac{1}{2}+\bb^T\bb$, $\lambda_3=\frac{1}{2}$ and the eigenvalues of $B$ are $\lambda_1=-\frac{1}{2}-\cb^T\cb$, $\lambda'_2=-\frac{1}{2}-\db^T\db$ and $\lambda'_3=-\frac{1}{2}$. This follows from $\ab\ab^T+\bb\bb^T$ and $\cb\cb^T+\db\db^T$ being symmetric rank 2 matrices and the pairwise orthogonality between the vectors. It is therefore evident from the fADI error formula
\begin{equation*}   \min_{\substack{\alpha_i\in\Cbb\\\beta_i\in\Cbb}}\max_{\substack{\lambda_i\in\operatorname{eig}(A)\\\lambda'_i\in\operatorname{eig}(B)}}\prod_{i=1}^n\bigg|\frac{(\alpha_i-\lambda_i)(\alpha_i-\lambda'_i)}{(\alpha_i-\lambda'_i)(\beta_i-\lambda_i)}\bigg|,
\end{equation*}
that choosing the shift parameters as $\alpha_j=\lambda_j$ and $\beta_j=\lambda'_j$ makes the iteration converge in three iterations. We can therefore represent the solution using the fADI solution formula as

$$X=\Tc_\epsilon^{sum}\bigg(\sum_{j=1}^3Z_i(\beta_i-\alpha_i)I_{r\times r}Y_i\bigg) $$
and what remains is to determine the factors $Z_i,Y_i$ for $i=1,2,3$. Before deriving the explicit forms we need the following helpful calculation. The inverse of an expression of the form $(\vb^T\vb+1)I_n+\ub\ub^T$ for arbitrary vectors $\ub,\vb\in\Rbb^n$ can be obtained using the Sherman-Morrison formula \cite{Petersen2008} as 
\begin{equation}
   [(1+\vb^T\vb)I_n+\ub\ub^T]^{-1}= \frac{1}{1+\vb^T\vb}\bigg[I_n-\frac{\ub\ub^T}{1+\vb^T\vb+\ub^T\ub}\bigg].
   \label{eq:helpresultfADI}
\end{equation}
We now derive the explicit forms of the matrices $Z_i$, omitting the calculations for $Y_i$ as they are constructed similarly. Recall from \cite{Benner2009} that $Z_1$ is given by $Z_1=(A-\beta_1 I_n)^{-1}U\sqrt{S}$. Straightforward calculations give

\begin{align*}
    Z_1&=((1+\cb^ T\cb)I_n+\ab\ab^T+\bb\bb^T)^{-1}U\sqrt{S}\\
    &=\bigg[C_1^{-1}-\frac{C^{-1}\bb\bb^TC^{-1}}{1+\bb^TC^{-1}\bb}\bigg]U\sqrt{S},
\end{align*}
where $C=(1+\cb^T\cb)I_n+\ab\ab^T$. Using \eqref{eq:helpresultfADI} together with the fact that $\ab^T\bb=0$ we then get the final form of $Z_1$ as
\begin{align*}
    Z_1&=\bigg[\frac{1}{1+\cb^T\cb}I_n-\frac{\ab\ab^T}{1+\ab^T\ab+\cb^T\cb}-\frac{\bb\bb^T}{(1+\cb^T\cb)(1+\cb^T\cb+\bb^T\bb}\bigg]U\sqrt{S}\\
    &=\begin{bmatrix}
        U & \ab & \bb
    \end{bmatrix}
        \begin{bmatrix}
            \frac{\sqrt{S}}{1+\cb^T\cb} & \gamb_a^T\sqrt{S} & \gamb_b^T\sqrt{S}
        \end{bmatrix}
\end{align*}
where we for notational convenience introduce the vectors
\begin{align*}
    &\gamb_a^T=-\frac{\ab^TU}{(1+\cb^T\cb)(1+\ab^T\ab+\cb^T\cb)},\\
    &\gamb_b=-\frac{\bb^TU}{(1+\cb^T\cb)(1+\bb^T\bb+\cb^T\cb)}.
\end{align*}
Then, $Z_2$ is obtained by

\begin{align*}
    Z_2&=Z_1+(\beta_2-\alpha_1)(A-\beta_2I_n)^{-1}Z_1\\
    &=[I_n-(1+\ab^T\ab+\db^T\db)((1+\db^T\db)I_n+\ab\ab^T+\db\db^T)^{-1}]Z_1\\
    &=\bigg[-\frac{\ab^T\ab}{1+\db^T\db}I_n+\frac{1}{1+\db^T\db}\ab\ab^T+\frac{1+\ab^T\ab+\db^T\db}{1+\bb^T\bb+\db^T\db}\frac{\bb\bb^T}{1+\db^T\db}\bigg]Z_1\\
    &=[\kz I+\kzz\ab\ab^T+\kzzz\bb\bb^T]Z_1
\end{align*}
where equation \eqref{eq:helpresultfADI} was used twice together with the orthogonality between $\ab$ and $\bb$ in the second step. Finally, $Z_3$ follows from a similar calculation

\begin{align*}
    Z_3&=Z_2+(\beta_3-\alpha_2)(A-\beta_3I_n)^{-1}Z_2\\
    &=[I_n-(1+\bb^T\bb)(I_n+\ab\ab^T+\bb\bb^T)^{-1}]Z_2\\
    &=\bigg[-\bb^TI_n+\frac{1+\bb^T\bb}{1+\ab^T\ab}\ab\ab^T+\bb\bb^T\bigg]Z_2\\
    &=\bigg[-\frac{\bb^T\bb}{1+\db^T\db}I_n-\frac{\bb^T\bb}{1+\db^T\db}\ab\ab^T-\frac{\ab^T\ab}{1+\db^T\db}\bb\bb^T\bigg]Z_1\\
    &=[\kzt I_n+\kzzt\ab\ab^T+\kzzzt\bb\bb^T]Z_1,
\end{align*}
and the lemma follows by analogous calculations for $Y_i$.

\section{Derivation for Problem~\ref{Lemma:LRAAweights} }
\label{proof:LRAA}
\begin{proof}
    We derive the system of equations by a straightforward calculation. Note that if $\gamb^{(k)}$ is the desired minimizer, it necessarily minimizes the internal sum of squared norms $f(\ub)=\sum_{l=1}\|D^k_l\Lambda(\ub)-F^k_l\|^2$. We can therefore view this as a multivariate function of the variables in $\ub$. To solve the minimization problem we calculate the components of $\nabla_\ub f(\ub)$ and equate them to zero. To make the calculations more transparent, we first rewrite $f(\ub)$ by expanding the inner product $D^k_l\Lambda(\ub)$ and writing out the norm in index notation 
    \begin{align*}
      f(\ub)=\sum_{l=1}^p\bigg \|\sum_{j=1}^m\Delta F^{k-j}_lu_j-F^k_l\bigg \|^2=\sum_{l=1}^p\bigg(\sum_{\alpha,\beta=1}^n\bigg[\sum_{j=1}^m\Delta F^{k-j}_l(\alpha,\beta)u_j-F^k_l(\alpha,\beta)\bigg]^2\bigg).
    \end{align*}
    Differentiating $f(\ub)$ with respect to a component $u_i$ we then get
    \begin{align*}
        \frac{1}{2}\frac{\partial f(\ub)}{\partial u_i}=&\sum_{l=1}^p\bigg(\sum_{\alpha,\beta=1}^n\bigg[\sum_{j=1}^m\Delta F^{k-j}_l(\alpha,\beta)u_j-F^k_l(\alpha,\beta)\bigg]\Delta F^{k-j}_l\frac{\partial u_j}{\partial u_i}\bigg)\\
        &=\sum_{l=1}^p\bigg(\sum_{\alpha,\beta=1}^n\bigg[\sum_{j=1}^m\Delta F^{k-j}_l(\alpha,\beta)\Delta F^{k-i}_l(\alpha,\beta)u_j-F^k_l(\alpha,\beta)\Delta F^{k-i}_l(\alpha,\beta)\bigg]\bigg)\\
        &=\sum_{l=1}^p\bigg(\sum_{j=1}^m\langle\Delta F^{k-i}_l,\Delta F^{k-j}_l\rangle u_j-\langle F^k_l,\Delta F^{k-i}_l\rangle\bigg)=0,
    \end{align*}
    where we in the first step used that $\frac{\partial u_j}{\partial u_i}$ becomes the Kronecker delta and in the second step introduced a low-rank Frobenius inner product $\langle\cdot,\cdot\rangle$. Realizing that the first term in the final parenthesis corresponds to the $i$'th row of a matrix-vector product $A\gamb$ with the matrix $A$ given by  $A_{ij}=\sum_{l=1}^p\langle \Delta F^{k-i}_l,\Delta F^{k-j}_l\rangle$ and the right-hand side the $i$'th component of a vector $\bb$ given by $ b_i=\sum_{l=1}^p\langle F^k_l,\Delta F^{k-i}_l\rangle$ the desired solution follows. Performing one more differentiation makes it clear that $\gamb$ yields a minimum as $\frac{\partial^2f(\ub)}{\partial u_{i}^2}>0$ for all components $u_i$.
\end{proof}

\end{document}